\documentclass[12pt]{amsart}
\usepackage{a4wide,amsfonts,amssymb,stmaryrd,
amscd,amsmath,latexsym,amsbsy}
\numberwithin{equation}{section}
\theoremstyle{plain}
\newtheorem{thm}{Theorem}[section]

\newtheorem{prop}[thm]{Proposition}

\newtheorem{defi}[thm]{Definition}
\newtheorem{lem}[thm]{Lemma}
\newtheorem{cor}[thm]{Corollary}

\theoremstyle{remark}
\newtheorem{rema}[thm]{Remark}
\newtheorem{eg1}[thm]{{Example}}
\title[Vertex-IRF transformations 
and harmonic analysis]{Vertex-IRF transformations, 
dynamical quantum groups and harmonic analysis} 
\author{Jasper V. Stokman}
\address{J.V. Stokman,
KdV Institute for Mathematics, Universiteit van Amsterdam,
Plantage Muidergracht 24, 1018 TV Amsterdam, The Netherlands.}
\email{jstokman@science.uva.nl}
\begin{document}

\begin{abstract}
It is shown that a dynamical quantum group arising from a 
vertex-IRF transformation
has a second realization with untwisted dynamical multiplication
but nontrivial bigrading. Applied to the $\hbox{SL}(2;\mathbb{C})$ dynamical
quantum group, the second realization is naturally described in terms
of Koornwinder's twisted primitive elements. This 
leads to an intrinsic explanation why harmonic analysis on the ``classical''
$\hbox{SL}(2;\mathbb{C})$ quantum group
with respect to twisted primitive elements, as initiated
by Koornwinder, is the same as harmonic analysis 
on the $\hbox{SL}(2;\mathbb{C})$
dynamical quantum group.
\end{abstract}

\maketitle
\begin{center}
{\it Dedicated to Tom Koornwinder on the occasion of his 60th birthday}
\end{center}


\section{Introduction}

In the remarkable paper \cite{Koo1}, Koornwinder
obtained in 1993 a subclass of the most general, classical family of
basic hypergeometric orthogonal polynomials, known as Askey-Wilson \cite{AW}
polynomials, as spherical functions on the $\hbox{SL}(2;\mathbb{C})$
quantum group. Koornwinder's
results have been refined (see e.g. \cite{NM}, \cite{Koe}, \cite{Koe2}) 
and have been succesfully extended to higher rank (see e.g.
\cite{N}, \cite{NDS}), leading to the interpretation
of Macdonald-Koornwinder polynomials \cite{Mac}, \cite{K} as spherical
functions on higher rank quantum symmetric spaces.
In this paper we provide a natural reinterpretation for some of
Koornwinder's \cite{Koo1} techniques and results in terms of
the {\it trigonometric} $\hbox{SL}(2;\mathbb{C})$ 
{\it dynamical quantum group}.

In 1998 and 1999, Etingof and Varchenko \cite{EV1}, \cite{EV2}
introduced the notion of a dynamical quantum group. A dynamical 
quantum group is a Hopf algebroid obtained by twisting
a Hopf algebra by a {\it dynamical twist}. Dynamical 
twists occur naturally in conformal field theory as universal fusion 
matrices, see e.g. \cite{EV3}. 
In \S 3 we give the construction of a 
dynamical quantum group starting 
from a five-tuple $(\mathcal{U},\mathcal{A},
\langle \cdot,\cdot\rangle,T,J(\lambda))$ 
with $\langle \cdot,\cdot\rangle$
a Hopf algebra pairing between two 
Hopf algebras $\mathcal{U}$ and $\mathcal{A}$,
$T$ a finite abelian subgroup of the group-like 
elements of $\mathcal{U}$ and $J(\lambda)$ a dynamical 
twist for $\mathcal{U}$ with respect to $T$. The dynamical
quantum group construction in \cite{EV2} is the special case that
$\mathcal{A}$ is the space of matrix coefficients of a suitable
tensor category of $\mathcal{U}$-representations.

We are specifically interested in dynamical quantum groups with
associated dynamical twists arising from 
{\it vertex-IRF transformations} \cite{EN2}.
In statistical mechanics vertex-IRF transformations arise 
as gauge transformations
relating vertex models to Interaction-Round-a-Face (IRF) models. 
The notions of 
vertex-IRF transformations and
dynamical twists are recalled in \S 2.
We show in \S 3 that the dynamical quantum group arising from a vertex-IRF
transformation has, 
besides its usual realization with twisted dynamical
multiplication and trivial bigrading, 
a second realization with untwisted dynamical multiplication
and twisted bigrading. 

The trigonometric $\hbox{SL}(2;\mathbb{C})$ dynamical quantum group
is one of the few known nontrivial examples of a dynamical quantum group
arising
from a vertex-IRF transformation \cite{B}, \cite{BBB}. 
We show in \S 4 that the alternative realization of the trigonometric 
$\hbox{SL}(2;\mathbb{C})$
dynamical quantum group, with untwisted dynamical multiplication and twisted
bigrading, 
can be naturally formulated in terms of the eigenspace decomposition  
of the $\hbox{SL}(2;\mathbb{C})$ quantum group with respect to 
Koorwinder's \cite{Koo1} twisted primitive elements.
The key fact needed here is the observation of
Rosengren \cite{R0} that the vertex-IRF transformation conjugates
a standard Cartan type element to Koornwinder's  
twisted primitive elements.

As a consequence, we obtain an intrinsic link between
Koelink's and Rosengren's \cite{KR} harmonic analysis on the 
trigonometric $\hbox{SL}(2;\mathbb{C})$
dynamical quantum group and
harmonic analysis on 
the $\hbox{SL}(2;\mathbb{C})$ quantum group as initiated by
Koornwinder \cite{Koo1} and extended 
by Noumi, Mimachi \cite{NM} and Koelink \cite{Koe}.
Such link was predicted by Koelink and Rosengren
in the introduction of \cite{KR} after a direct comparison of
their harmonic analytic results 
to the corresponding results in \cite{Koo1},
\cite{NM} and \cite{Koe}. 
We explore this intrinsic link in some detail by deriving 
several harmonic
analytic results 
from \cite{Koo1}, \cite{NM} and \cite{Koe} as a direct consequence
of the corresponding results on the $\hbox{SL}(2;\mathbb{C})$ 
dynamical quantum group from \cite{KR}.

{\it Acknowledgments:}  
A substantial part of the research was done 
while visiting the Massachusetts Institute of Technology, USA
in 2002. I thank MIT for the hospitality and
Pavel Etingof and Alexei Oblomkov for discussions. 
I am supported by the Royal Netherlands Academy
of Arts and Sciences (KNAW).


\section{Dynamical twists and vertex-IRF transformations}

This section follows closely Etingof and Nikshych
\cite[\S 2.1]{EN2}, see also \cite[\S 4.1]{EN1}.

Let $(\mathcal{U},m_{\mathcal{U}},1_{\mathcal{U}},
\Delta_{\mathcal{U}},\epsilon_{\mathcal{U}},S_{\mathcal{U}})$ 
be a Hopf-algebra over $\mathbb{C}$
with multiplication $m_{\mathcal{U}}$, unit $1_{\mathcal{U}}$, 
comultiplication $\Delta_{\mathcal{U}}$,
counit $\epsilon_{\mathcal{U}}$ and antipode $S_{\mathcal{U}}$. 
For $u\in \mathcal{U}^{\otimes k}$ with $k\leq n$ and for an 
ordered $k$-tuple $\{i_1,\ldots,i_k\}\subset \{1,\ldots,n\}$ we
write $u_{i_1i_2\cdots i_k}$ for the element of $\mathcal{U}^{\otimes n}$
obtained by mapping the $j$th tensor component of $\mathcal{U}^{\otimes k}$
to the $i_j$th tensor component of $\mathcal{U}^{\otimes n}$,
where $\otimes=\otimes_{\mathbb{C}}$ denotes the usual tensor product 
over $\mathbb{C}$.
For instance, $(u\otimes v)_{31}=v\otimes 1\otimes u$ 
when regarded as element in $\mathcal{U}^{\otimes 3}$.

For $j\in \{1,\ldots,n\}$ we define algebra homomorphisms 
$\Delta_{\mathcal{U}}{}_j: \mathcal{U}^{\otimes n}\rightarrow 
\mathcal{U}^{\otimes (n+1)}$ and
$\epsilon_{\mathcal{U}}{}_j: {\mathcal{U}}^{\otimes n}\rightarrow 
\mathcal{U}^{\otimes (n-1)}$ by letting
$\Delta_{\mathcal{U}}$ and $\epsilon_{\mathcal{U}}$ 
act on the $j$th tensor component.
Similarly we define the map $S_{\mathcal{U}j}: 
\mathcal{U}^{\otimes n}\rightarrow 
\mathcal{U}^{\otimes n}$ by letting the antipode
$S_{\mathcal{U}}$ act on the $j$th tensor component.
With these notations, the Hopf-algebra identities are 
\begin{equation}\label{Hopfaxioms}
\begin{split}
&\Delta_{\mathcal{U}1}\circ\Delta_{\mathcal{U}}=
\Delta_{\mathcal{U}2}\circ\Delta_{\mathcal{U}},\\
&\epsilon_{\mathcal{U}1}\circ\Delta_{\mathcal{U}}=
\hbox{Id}_{\mathcal{U}}=
\epsilon_{\mathcal{U}2}\circ\Delta_{\mathcal{U}},\\
&m_{\mathcal{U}}\circ S_{\mathcal{U}1}\circ \Delta_{\mathcal{U}}(\cdot)=
\epsilon_{\mathcal{U}}(\cdot)1_{\mathcal{U}}=
m_{\mathcal{U}}\circ S_{\mathcal{U}2}\circ \Delta_{\mathcal{U}}(\cdot).
\end{split}
\end{equation}
We define the iterated comultiplication
$\Delta_{\mathcal{U}}^{(n-1)}: \mathcal{U}\rightarrow
\mathcal{U}^{\otimes n}$ inductively by
$\Delta_{\mathcal{U}}^{(0)}=\hbox{Id}_{\mathcal{U}}$ for $n=1$ and
$\Delta_{\mathcal{U}}^{(n-1)}=\Delta_{\mathcal{U}1}\circ 
\Delta_{\mathcal{U}}^{(n-2)}$
for $n\in\mathbb{Z}_{>1}$.

For a unital associative $\mathbb{C}$-algebra $L$ we denote
$\mathcal{U}_L=L\otimes \mathcal{U}$ for the Hopf algebra over $L$
obtained by extending the Hopf-algebra maps $L$-linearly.
We keep the notation $1_{\mathcal{U}}$ for the unit of $\mathcal{U}_L$
and $m_{\mathcal{U}}$, $\Delta_{\mathcal{U}}$, $\epsilon_{\mathcal{U}}$ and
$S_{\mathcal{U}}$ for the $L$-linear 
extended Hopf algebra maps of $\mathcal{U}_L$.

Let $\mathcal{U}^\times$ be the multiplicative group of 
invertible elements in $\mathcal{U}$. The {\it group-like} elements
\[
G(\mathcal{U})=\{u\in \mathcal{U}\, | \, 
\Delta_{\mathcal{U}}(u)=u\otimes u,\,\,\, \epsilon_{\mathcal{U}}(u)=1\}
\]
is a subgroup of $\mathcal{U}^\times$ since 
$S_{\mathcal{U}}(u)$ is the inverse of a group-like element $u$.
Let $T\subseteq G(\mathcal{U})$ be a finite abelian subgroup and
write $\widehat{T}$ for the character group of $T$. 
The value of a character $\alpha\in\widehat{T}$ at $t\in T$ is denoted
by $t^\alpha$. 

We suppose that $\mathcal{U}$ is $\hbox{ad}(T)$-semisimple, where
\[\hbox{ad}(u)v=\sum u_{(1)}vS_{\mathcal{U}}(u_{(2)}),\qquad
u,v\in \mathcal{U}
\]
is the analogue of the adjoint action.
Consequently, $\mathcal{U}$ is $\widehat{T}$-graded,
\[\mathcal{U}=\bigoplus_{\alpha\in\widehat{T}}
\mathcal{U}[\alpha]
\]
with $\mathcal{U}[\alpha]$ the elements $u\in \mathcal{U}$ satisfying
$\hbox{ad}(t)u=t^\alpha u$ for all $t\in T$. Since $T$ is abelian we
have $T\subseteq \mathcal{U}[0]$. Furthermore, for $\alpha,\beta
\in\widehat{T}$ with $\alpha\not=0$,
\[\epsilon_{\mathcal{U}}(\mathcal{U}[\alpha])=0,
\qquad \Delta_{\mathcal{U}}(\mathcal{U}[\beta])\subseteq
\bigoplus_{\gamma\in \widehat{T}}\mathcal{U}[\gamma]\otimes 
\mathcal{U}[\beta-\gamma],\qquad
S_{\mathcal{U}}\bigl(\mathcal{U}[\beta]\bigr)\subseteq
\mathcal{U}[\beta].
\]

The primitive idempotents of $T$ are 
\[\pi_{\alpha}=
\frac{1}{\#T}\sum_{t\in T}t^{-\alpha}\,t\in\mathcal{U}\qquad
(\alpha\in\widehat{T}).
\]
Their basic properties are
\[\sum_{\gamma\in\widehat{T}}\pi_\gamma=1_{\mathcal{U}},\qquad
t\pi_\alpha=t^\alpha\pi_\alpha=\pi_\alpha t,\qquad
u\pi_\alpha=\pi_{\alpha+\beta}u
\]
for $t\in T$, $u\in\mathcal{U}[\beta]$ and $\alpha,\beta\in\widehat{T}$,
and
\begin{equation*}
\begin{split}
\pi_\alpha\pi_\beta=
\delta_{\alpha,\beta}\pi_\alpha,\qquad
&\Delta_{\mathcal{U}}(\pi_\alpha)=
\sum_{\beta\in\widehat{T}}\pi_\beta\otimes
\pi_{\beta-\alpha},\\
\epsilon_{\mathcal{U}}(\pi_\alpha)=\delta_{\alpha,0},\qquad
&S_{\mathcal{U}}(\pi_\alpha)=\pi_{-\alpha}
\end{split}
\end{equation*}
for $\alpha,\beta\in \widehat{T}$ with $\delta_{\alpha,\beta}$ the
Kronecker delta function on $\widehat{T}\times\widehat{T}$.

Let $F$ be the $\mathbb{C}$-algebra of complex valued 
functions on $\widehat{T}$
and set $K=F\otimes F$. An element $f\in K$ is denoted
by $f(\lambda,\mu)=\sum_jf_j(\lambda)\otimes f_j^\prime(\mu)$, with 
$f_j,f_j^\prime\in F$. The variables $\lambda,\mu\in\widehat{T}$ 
thus indicate the 
dependence in the first and second tensor
component of $K=F\otimes F$, respectively. In particular, we have
the subalgebras $F_1$ and $F_2$ of $K$ consisting of elements
$f(\lambda)=f\otimes 1$ and $g(\mu)=1\otimes g$, respectively. 
Similar notations will be used
for elements in $\mathcal{U}_K$,  
$\mathcal{U}_{F_1}$ and $\mathcal{U}_{F_2}$. 

For $u(\lambda,\mu)\in \mathcal{U}^{\otimes n}_K$ 
and $j\in \{1,\ldots,n\}$ we denote
\begin{equation*}
\begin{split}
u(\lambda\pm h^{(j)},\mu)&=\sum_{\alpha\in\widehat{T}}
u(\lambda\pm\alpha,\mu)\pi_{{\alpha}j}
\in\mathcal{U}^{\otimes n}_{K},\\
u\bigl(\lambda,\mu\pm h^{(j)}\bigr)&=
\sum_{\alpha\in\widehat{T}}
u(\lambda,\mu\pm\alpha)\pi_{{\alpha}j}\in \mathcal{U}^{\otimes n}_{K}.
\end{split}
\end{equation*}
We are now in a position to give the definition of a {\it dynamical twist},
see \cite[Def. 2.3]{EN2}.
\begin{defi}
{\bf (i)} We say that $u\in \mathcal{U}_K^{\otimes n}$ 
is of zero $T$-weight if 
\[\lbrack u,\Delta_{\mathcal{U}}^{(n-1)}(t)\rbrack=0,\qquad\forall\, 
t\in T.
\]
{\bf (ii)}
An invertible element 
$J(\lambda)\in \mathcal{U}_{F_1}^{\otimes 2}$
is called a {\it dynamical twist} with respect to $T$
if $J(\lambda)$ is of zero $T$-weight satisfying
\begin{equation}\label{dyncond}
\begin{split}
&\epsilon_{\mathcal{U}1}(J(\lambda))=1_{\mathcal{U}}=
\epsilon_{\mathcal{U}2}(J(\lambda)),\\
&\Delta_{\mathcal{U}1}(J(\lambda))J_{12}(\lambda+h^{(3)})=
\Delta_{\mathcal{U}2}(J(\lambda))J_{23}(\lambda).
\end{split}
\end{equation}
\end{defi}
Observe that $\mathcal{U}_K[0]$ is the 
set of zero $T$-weighted elements in
$\mathcal{U}_K$, and  
$\bigoplus_{\alpha\in\widehat{T}}K\otimes \mathcal{U}[\alpha]\otimes
\mathcal{U}[-\alpha]$ is the set of zero $T$-weighted elements in 
$\mathcal{U}^{\otimes 2}_K$.

\begin{eg1}\label{Jexample}
Suppose $\mathcal{U}$ is the 
Drinfeld-Jimbo quantized universal enveloping
algebra of a complex semisimple Lie algebra. 
A dynamical twist $J(\lambda)$
then naturally arises
as the universal fusion 
matrix of $\mathcal{U}$, see \cite{EV3} and \cite{KT}.
\end{eg1}

\begin{defi}\label{gauge}
Let $x(\lambda)\in\mathcal{U}_{F_1}$ be invertible
and suppose that $\epsilon_{\mathcal{U}}(x(\lambda))=1$.
If $x(\lambda)$ is of zero $T$-weight then $x(\lambda)$ is called
a gauge transformation with respect to $T$.
\end{defi}
Gauge transformations can be used to ``gauge'' dynamical twists,
see \cite[\S 2.1]{EN2}.
\begin{lem}
Let $J(\lambda)\in \mathcal{U}^{\otimes 2}_{F_1}$ be a dynamical
twist and $x(\lambda)\in \mathcal{U}_{F_1}$ be a gauge 
transformation with respect to $T$.
Then
\begin{equation}\label{gaugeJ}
J_x(\lambda)=\Delta_{\mathcal{U}}(x(\lambda))J(\lambda)
x_2^{-1}(\lambda)x_1^{-1}(\lambda+h^{(2)})
\end{equation}
is a dynamical twist with respect to $T$.
\end{lem}
On the other hand, 
a dynamical twist ``almost'' gives rise to a gauge transformation,
cf. \cite[\S 4.2]{EV2}, \cite[\S 4.3]{EN1} and \cite[Lem. 2.12]{EV3}.
\begin{lem}\label{KQ}
If $J(\lambda)\in\mathcal{U}_{F_1}^{\otimes 2}$ is a dynamical
twist with respect to $T$, then
\[K^J(\lambda)=m_{\mathcal{U}}\bigl(S_{\mathcal{U}1}(J(\lambda))\bigr),\qquad
Q^J(\lambda)=m_{\mathcal{U}}\bigl(S_{\mathcal{U}2}(J^{-1}(\lambda))\bigr)
\]
are of zero $T$-weight and satisfy
\begin{equation}\label{inverse}
Q^J(\lambda+h)K^J(\lambda)=1_{\mathcal{U}},
\end{equation}
where $Q^J(\lambda\pm h)=
\sum_{\alpha\in\widehat{T}} Q^J(\lambda\pm\alpha)\pi_\alpha$.
\end{lem}
\begin{proof}
The only nontrivial part of 
Lemma \ref{KQ} is \eqref{inverse}, which follows from
applying $m_{\mathcal{U}}\circ \bigl(m_{\mathcal{U}}\otimes 
\hbox{Id}_{\mathcal{U}}\bigr)
\circ S_{\mathcal{U}2}$ to the reformulation 
\[
\Delta_{\mathcal{U}1}(J(\lambda))=\Delta_{\mathcal{U}2}(J(\lambda))
J_{23}(\lambda)J_{12}^{-1}(\lambda+h^{(3)})
\]
of the second equality of \eqref{dyncond}, 
and by using \eqref{Hopfaxioms} and
the first equality of \eqref{dyncond}.
\end{proof}

If $\mathcal{U}_{F_1}$ does not have zero divisors, then it follows 
from Lemma
\ref{KQ} that 
$K^J(\lambda)$ and $Q^J(\lambda)$ are gauge transformations 
with respect to $T$, and that $(K^J)^{-1}(\lambda)=Q^J(\lambda+h)$.

A vertex-IRF transformation is the following generalization 
of a gauge transformation, see \cite[Def. 2.6]{EN2}.
\begin{defi}
Let $x(\lambda)\in\mathcal{U}_{F_1}$ be 
invertible and suppose that 
$\epsilon_{\mathcal{U}}(x(\lambda))=1$. We say that $x(\lambda)$
is a vertex-IRF transformation with respect to $T$ if
\begin{equation}\label{Jx}
j_x(\lambda)=
\Delta_{\mathcal{U}}(x(\lambda))x_2^{-1}(\lambda)x_1^{-1}(\lambda+h^{(2)})
\in\mathcal{U}_{F_1}^{\otimes 2}
\end{equation}
is of zero $T$-weight.
\end{defi}
A vertex-IRF transformation gives rise to a dynamical twist,
see \cite[Prop. 2.5]{EN2}.
\begin{prop}
If $x(\lambda)\in\mathcal{U}_{F_1}$ is a vertex-IRF transformation with
respect to $T$ then
$j_x(\lambda)\in\mathcal{U}_{F_1}^{\otimes 2}$ is a dynamical
twist with respect to $T$.
\end{prop}
Key examples of vertex-IRF transformations arise 
in the theory of exactly solvable lattice models 
as gauge transformations relating vertex models to 
Interaction-Round-a-Face (IRF) models, see e.g. \cite{F} and references
therein. 

We end this section by considering the special case of quasi-triangular
Hopf-algebra's $\mathcal{U}$. In that case the 
Hopf algebra $\mathcal{U}$ has a universal $R$-matrix 
$\mathcal{R}\in\mathcal{U}\otimes \mathcal{U}$, which
is invertible and satisfies 
\[\Delta_{\mathcal{U}}^{op}(\cdot)=
\mathcal{R}\Delta_{\mathcal{U}}(\cdot)\mathcal{R}^{-1},\qquad
\Delta_{\mathcal{U}1}(\mathcal{R})=
\mathcal{R}_{13}\mathcal{R}_{23},\qquad
\Delta_{\mathcal{U}2}(\mathcal{R})= \mathcal{R}_{13}\mathcal{R}_{12},
\]
with $\Delta_{\mathcal{U}}^{op}(u)=
\bigl(\Delta_{\mathcal{U}}(u)\bigr)_{21}$ the opposite
comultiplication. 
In particular, $\mathcal{R}$ satisfies 
the quantum Yang-Baxter equation (QYBE)
\[\mathcal{R}_{12}\mathcal{R}_{13}\mathcal{R}_{23}=
\mathcal{R}_{23}\mathcal{R}_{13}\mathcal{R}_{12}
\]
in $\mathcal{U}^{\otimes 3}$. 
Gauging the universal $R$-matrix by a dynamical twist leads to a solution
of a dynamical version of the quantum
Yang-Baxter equation.
The precise result is as follows, see \cite[\S 3]{BBB}.
\begin{prop}\label{JtwistR}
Let $J(\lambda)\in\mathcal{U}^{\otimes 2}_{F_1}$ 
be a dynamical twist with respect to $T$.
Then 
\begin{equation}\label{Rdyn}
\mathcal{R}^J(\lambda)=J^{-1}(\lambda)\mathcal{R}_{21}J_{21}(\lambda)\in
\mathcal{U}^{\otimes 2}_{F_1}
\end{equation}
satisfies
\begin{equation}\label{DQYBE}
\mathcal{R}_{12}^J(\lambda+h^{(3)})\mathcal{R}_{13}^J(\lambda)
\mathcal{R}_{23}^J(\lambda+h^{(1)})=
\mathcal{R}_{23}^J(\lambda)\mathcal{R}_{13}^J(\lambda+h^{(2)})
\mathcal{R}_{12}^J(\lambda).
\end{equation}
\end{prop}
As Koornwinder emphasizes in \cite{Koo2}, 
the proof of Proposition \ref{JtwistR}
is a direct consequence of the identities
\begin{equation}\label{intermediate}
\begin{split}
J_{12}^{-1}(\lambda+h^{(3)})
\Delta_{\mathcal{U}1}(\mathcal{R}^J(\lambda))J_{12}(\lambda)&
=\mathcal{R}_{23}^J(\lambda)\mathcal{R}_{13}^J(\lambda+h^{(2)}),\\
J_{23}^{-1}(\lambda)
\Delta_{\mathcal{U}2}(\mathcal{R}^J(\lambda))J_{23}(\lambda+h^{(1)})&
=\mathcal{R}_{12}^J(\lambda+h^{(3)})\mathcal{R}_{13}^J(\lambda)
\end{split}
\end{equation}
which are dynamical analogues of 
$\Delta_{\mathcal{U}1}(\mathcal{R}_{21})=\mathcal{R}_{32}\mathcal{R}_{31}$
and $\Delta_{\mathcal{U}2}(\mathcal{R}_{21})=
\mathcal{R}_{21}\mathcal{R}_{31}$ 
respectively. 

\begin{defi}
Equation \eqref{DQYBE} is called the dynamical
quantum Yang-Baxter equation \textup{(}DQYBE\textup{)}. A solution 
$\mathcal{R}(\lambda)\in \mathcal{U}_{F_1}^{\otimes 2}$ is called
a dynamical universal $R$-matrix.
\end{defi}
The DQYBE is also known as the Gervais-Neveu-Felder equation,
and is closely related to Baxter's star-triangle relation, 
see e.g. \cite{GN} and \cite{F}. 
\begin{eg1}\label{Rexample}
In the setup of example \ref{Jexample}, $\mathcal{U}$ is quasi-triangular,
and the universal $R$-matrix $\mathcal{R}^J(\lambda)$ associated to the
universal fusion matrix $J(\lambda)$
is the universal 
exchange matrix of $\mathcal{U}$, see \cite{EV3} and \cite{KT}.
\end{eg1}

The complete integrability of vertex models and IRF models are governed
by the quantum Yang-Baxter equation and Baxter's star-triangle identity,
respectively. For quasi-triangular Hopf-algebras $\mathcal{U}$
with vertex-IRF transformation $x(\lambda)$, Proposition \ref{JtwistR}
shows that these basic integrability conditions
are interrelated by twisting 
the corresponding universal $R$-matrix $\mathcal{R}$
of $\mathcal{U}$ with the dynamical twist $j_x(\lambda)$. Note that 
the corresponding universal dynamical $R$-matrix 
$\mathcal{R}^{j_x}(\lambda)=
j_x^{-1}(\lambda)\mathcal{R}_{21}(j_x)_{21}(\lambda)$
can alternatively be written as
\begin{equation}\label{Rdynx}
\mathcal{R}^{j_x}(\lambda)=x_1(\lambda+h^{(2)})x_2(\lambda)\mathcal{R}_{21}
x_1^{-1}(\lambda)x_2^{-1}(\lambda+h^{(1)}),
\end{equation}
see \cite[Cor. 2.11]{EN2}.


\section{Dynamical quantum groups}

Etingof and Varchenko \cite{EV1}, \cite{EV2}
gave a general construction of a Hopf algebroid starting from
a given nondegenerate, polarized Hopf algebra $\mathcal{U}$
and a suitable tensor category of $\mathcal{U}$-representations. 
The notion of a Hopf algebroid 
is closely related to weak Hopf algebras and quantum groupoids, 
see \cite{EN2}. 
The constructed Hopf algebroids are modeled on 
the space of matrix coefficients
of the $\mathcal{U}$-representations from the given tensor category. 
The Hopf algebroid structures
are governed by the fusion matrices of $\mathcal{U}$, or, when $\mathcal{U}$
is quasi-triangular, by the exchange matrices for $\mathcal{U}$. These
Hopf algebroids are called 
{\it exchange dynamical quantum groups}, or simply {\it dynamical quantum
groups}. A different, but closely related construction was given in
\cite[\S 4.3]{EN1}.

In this section we give the construction of dynamical quantum groups 
in a slightly different setup. The input data is a five-tuple
$(\mathcal{U},\mathcal{A},\langle \cdot,\cdot\rangle,
T,J(\lambda))$
with $\mathcal{U}$ and $\mathcal{A}$ Hopf 
algebras, 
with $\langle \cdot,\cdot\rangle$ a Hopf pairing between $\mathcal{U}$
and $\mathcal{A}$, with $T\subset G(\mathcal{U})$
a finite abelian subgroup such that $\mathcal{U}$ is $\hbox{ad}(T)$-semisimple
and such that $\mathcal{A}$ is $(T-T)$-semisimple with respect to the 
left and right regular $\mathcal{U}$-action on $\mathcal{A}$, and with
$J(\lambda)\in\mathcal{U}^{\otimes 2}_{F_1}$ 
a dynamical twist for $\mathcal{U}$ with respect to $T$.

The resulting Hopf algebroid $\mathcal{A}^J$ is now
modeled on the space $\mathcal{A}_K=K\otimes \mathcal{A}$. 
We study $\mathcal{A}^J$ 
in more detail when the dynamical twist $J(\lambda)$ arises 
from a vertex-IRF transformation.


\subsection{Hopf algebroids}

This subsection follows closely \cite[\S 3.1]{EV2},
see also \cite[\S 2.1]{KR}, \cite[\S 2.2]{EN1} and references therein.

Let $\widehat{T}$ be the character group of some finite abelian group $T$,
with group operation written 
additively. Let $F$ be the unital $\mathbb{C}$-algebra
of complex valued functions on $\widehat{T}$.
We denote $1_F$ for the unit of $F$.
Translation over  $\alpha\in \widehat{T}$, 
\[\bigl(T_\alpha f\bigr)(\lambda)=f(\lambda+\alpha),\qquad f\in F,
\]
defines an automorphism of $F$. 
\begin{defi}
A $T$-algebra is a complex associative algebra $A$ with unit 
$1_A$
which is $\widehat{T}$-bigraded,
\[A=\bigoplus_{\alpha,\beta\in\widehat{T}}A_{\alpha\beta},
\]
and which is endowed with two unital algebra embeddings $\mu_l=\mu_l^A,
\mu_r=\mu_r^A: F\rightarrow
A_{00}$ satisfying
\begin{equation*}
\begin{split}
\mu_l(f)\circ\mu_r(g)&=\mu_r(g)\circ\mu_l(f),\\
\mu_l(f)\circ a&=a\circ\mu_l(T_{\alpha}f),\\
\mu_r(f)\circ a&=a\circ\mu_r(T_{\beta}f)
\end{split}
\end{equation*}
for $f,g\in F$ and $a\in A_{\alpha\beta}$, where $\circ$ denotes
the multiplication of the algebra $A$.
The algebra embeddings $\mu_l$ and $\mu_r$ are called the left and right
moment maps of $A$, respectively. 
\end{defi}
A morphism $\phi:A\rightarrow B$ between $T$-algebras $A$ and $B$
is an algebra homomorphism satisfying
\[\phi(A_{\alpha\beta})\subseteq B_{\alpha\beta},
\quad \phi\bigl(\mu_l^A(f)\bigr)=\mu_l^B(f),\quad
\phi\bigl(\mu_r^A(f)\bigr)=\mu_r^B(f)
\]
for $\alpha,\beta\in\widehat{T}$ and $f\in F$.

\begin{eg1}
The formal $|\widehat{T}|$-dimensional $F$-vector space 
$I=\bigoplus_{\alpha\in\widehat{T}}FT_\alpha$ is a unital, associative
$\mathbb{C}$-algebra with multiplication
$(fT_\alpha)\circ (gT_\beta)=\bigl(f(T_{\alpha}g)\bigr)T_{\alpha+\beta}$ and
unit $T_0$. It naturally acts on $F$ as difference operators on $\widehat{T}$
with coefficients from $F$. The algebra $I$ is a $T$-algebra 
with $(\alpha,\beta)$-bigraded piece 
\begin{equation*}
I_{\alpha\beta}=
\begin{cases}
0\qquad &\hbox{ when } \alpha\not=\beta,\\
FT_{-\alpha}\qquad &\hbox{ when } \alpha=\beta
\end{cases}
\end{equation*}
and moment maps $\mu_l(f)=\mu_r(f)=fT_0$.
\end{eg1}

Any $T$-algebra $A$ has the structure of a
$(F-F)$-bimodule,
\[ f\star a\star g=\mu_l(f)\circ\mu_r(g)\circ a,\qquad f,g\in F,\,\,a\in A.
\]
The bigraded pieces $A_{\alpha\beta}\subseteq A$ are $(F-F)$-submodules. 
We define the tensor product of two $T$-algebras $A$ and $B$
by
\[A\widetilde{\otimes} B=\bigoplus_{\alpha,\beta\in\widehat{T}}
\bigl(A\widetilde{\otimes} B\bigr)_{\alpha\beta},\qquad
\bigl(A\widetilde{\otimes} B\bigr)_{\alpha\beta}=
\bigoplus_{\gamma\in\widehat{T}}A_{\alpha\gamma}\otimes_FB_{\gamma\beta}.
\]
The balancing condition for the tensor product thus is 
$(\mu_r(f)\circ a)\otimes_Fb=
f\otimes_F(\mu_l(f)\circ b)$ for $f\in F$, $a\in A_{\alpha\gamma}$
and $b\in A_{\gamma\beta}$. Then $A\widetilde{\otimes}B$
becomes a $T$-algebra with multiplication
\[(a\otimes_Fb)\circ (a^\prime\otimes_Fb^\prime)=a\circ a^\prime
\otimes_Fb\circ b^\prime,
\]
with unit $1_A\otimes_F 1_B$, with $(\alpha,\beta)$-bigraded piece
$\bigl(A\widetilde{\otimes} B\bigr)_{\alpha\beta}$ and with moment maps
\[\mu_l(f)=\mu_l^A(f)\otimes_F1_B,\qquad \mu_r(f)=1_A\otimes_F\mu_r^B(f).
\]
We define for two morphisms $\phi: A\rightarrow A^\prime$ and 
$\psi: B\rightarrow B^\prime$ of $T$-algebras
a morphism  $\phi\widetilde{\otimes}\psi: 
A\widetilde{\otimes} A^\prime\rightarrow B\widetilde{\otimes} B^\prime$
by the usual formula
\[\bigl(\phi\widetilde{\otimes}\psi\bigr)(a\otimes_Fb)=
\phi(a)\otimes_F\psi(b).
\]
It is now straightforward to check that the category of $T$-algebras 
is a tensor category with tensor product
$\widetilde{\otimes}$, unit object $I$, the obvious associativity constraint
and unit constraints $l_A: I\widetilde{\otimes} A\rightarrow A$,
$r_A: A\widetilde{\otimes} I\rightarrow A$ given by
\[
l_A(fT_{-\alpha}\otimes_F a)=\mu_l^A(f)\circ a,\qquad
r_A(a\otimes_F fT_{-\beta})=\mu_r^A(f)\circ a,
\qquad f\in F,\,\,a\in A_{\alpha\beta}.
\]
In the remainder of the paper we use
the unit constraints to identify the $T$-algebras
$I\widetilde{\otimes}A$ and $A\widetilde{\otimes}I$ with $A$.
\begin{defi}
A $T$-bialgebroid is a $T$-algebra $A$
equipped with two morphisms $\Delta: A\rightarrow A\widetilde{\otimes} A$
and $\epsilon: A\rightarrow I$ satisfying the familiar coalgebra axioms
\[\bigl(\Delta\widetilde{\otimes}\hbox{Id}_A\bigr)\Delta=
\bigl(\hbox{Id}_A\widetilde{\otimes}\Delta\bigr)\Delta,\qquad
\bigl(\epsilon\widetilde{\otimes}\hbox{Id}_A\bigr)\Delta=
\hbox{Id}_A=\bigl(\hbox{Id}_A\widetilde{\otimes}\epsilon\bigr)\Delta.
\]
\end{defi}
The definition of a $T$-Hopf algebroid is a bit more subtle.
Suppose $A$ is a $T$-algebra and suppose that
$\phi: A\rightarrow A$
is a $\mathbb{C}$-linear map satisfying
\[\phi(\mu_r(f)\circ a)=\phi(a)\circ \mu_l(f),\qquad
\phi(a\circ \mu_l(f))=\mu_r(f)\circ \phi(a)
\]
for $a\in A$ and $f\in F$. Let
$\psi:A\rightarrow A$ be a morphism of $T$-algebras.
Then there exist unique $\mathbb{C}$-linear maps, 
denoted suggestively by
$m(\phi\otimes \psi),m(\psi\otimes\phi): A\widetilde{\otimes} 
A\rightarrow A$ with $m$ the multiplication map of $A$, such that
\[m\bigl(\phi\otimes \psi)(a\otimes_Fb)=\phi(a)\circ \psi(b),\qquad
m\bigl(\psi\otimes \phi)(a\otimes_Fb)=\psi(a)\circ \phi(b)
\]
for $a\in A_{\alpha\gamma}$ and $b\in A_{\gamma\beta}$.

For a difference operator $a\in I$ we denote $a1_F\in F$ for the
function obtained by applying $a$ to the constant function $1_F\in F$. 
In other words, $a1_F=\sum_{\alpha}a_\alpha(\lambda)$
when $a=\sum_{\alpha} a_\alpha(\lambda)T_{-\alpha}$. The following
definition of an antipode is from \cite[Def. 2.1]{KR}.

\begin{defi}
An antipode $S$ for a $T$-bialgebroid $A$ is a $\mathbb{C}$-linear
map $S: A\rightarrow A$ satisfying
\[S\bigl(\mu_r(f)\circ a)=S(a)\circ \mu_l(f),\qquad
S\bigl(a\circ \mu_l(f))=\mu_r(f)\circ S(a)
\]
for $f\in F$ and $a\in A$ and satisfying the antipode axioms
\[m(\hbox{Id}_A\otimes S)\bigl(\Delta(a)\bigr)=
\mu_l(\epsilon(a)1_F),\qquad
m\bigl(S\otimes \hbox{Id}_A\bigr)\bigl(\Delta(a)\bigr)=
\mu_r\bigl(T_\alpha(\epsilon(a)1_F)\bigr)
\]
for $a\in A_{\alpha\beta}$. The pair $(A,S)$ is called
a $T$-Hopf algebroid.
\end{defi}
\begin{defi}
A morphism $\phi: A\rightarrow B$ of 
$T$-Hopf algebroids $A$ and $B$ 
is a $T$-algebra morphism satisfying 
\[
\epsilon_B\bigl(\phi(a)\bigr)=\epsilon_A(a),\qquad
\Delta_B\bigl(\phi(a)\bigr)=
\bigl(\phi\widetilde{\otimes}\phi\bigr)\bigl(\Delta_A(a)\bigr),
\qquad
S_B\bigl(\phi(a)\bigr)=\phi\bigl(S_A(a)\bigr)
\]
for $a\in A$.
\end{defi}
 
For $T=\{1\}$ the trivial group, the definition of a
$T$-Hopf algebroid reduces to the familiar definition of
a Hopf algebra over $\mathbb{C}$. 


\subsection{The bigraded Hopf algebra}
Let $\mathcal{U}$ be a Hopf algebra over $\mathbb{C}$ as considered
in \S 2. Let $(\mathcal{A},m_{\mathcal{A}},1_{\mathcal{A}},
\Delta_{\mathcal{A}},\epsilon_{\mathcal{A}},S_{\mathcal{A}})$ be
a Hopf algebra over 
$\mathbb{C}$ and suppose that there exists an Hopf-algebra pairing 
$\langle \cdot,\cdot\rangle: \mathcal{U}\times
\mathcal{A}\rightarrow \mathbb{C}$ 
between $\mathcal{U}$ and $\mathcal{A}$, which we fix once and for all.
Then 
\begin{equation}\label{actiondot}
u\cdot a=\sum\langle u,a_{(2)}\rangle a_{(1)},
\qquad a\cdot u=\sum\langle u,a_{(1)}\rangle a_{(2)}
\end{equation}
for $u\in \mathcal{U}$ and $a\in\mathcal{A}$
defines a $(\mathcal{U}-\mathcal{U})$-bimodule structure on 
$\mathcal{A}$. We call \eqref{actiondot} the left and right regular action of
$\mathcal{U}$ on $\mathcal{A}$, respectively.

\begin{eg1}\label{example1}
Take $\mathcal{A}=\mathcal{U}^\star$ the Hopf-algebra dual of $\mathcal{U}$
with Hopf algebra pairing
\[\langle u,a\rangle=a(u),\qquad u\in\mathcal{U},\,\,a\in\mathcal{U}^\star.
\]
The associated $(\mathcal{U}-\mathcal{U})$-bimodule structure on 
$\mathcal{U}^\star$ is the regular action
\[\bigl(u\cdot a\cdot u^\prime\bigr)(v)=a(u^\prime vu),
\qquad a\in \mathcal{U}^\star,\quad u,u^\prime,v\in \mathcal{U}.
\]
\end{eg1}

Let $T\subseteq G(\mathcal{U})$ be a finite abelian subgroup such that
$\mathcal{U}$ is $\hbox{ad}(T)$-semisimple, cf. \S 2.
Recall that $K=F\otimes F$ with $F$ the algebra of complex valued functions
on $\widehat{T}$. We extend the Hopf-algebra maps of $\mathcal{U}$ and 
$\mathcal{A}$ $K$-linearly to arrive at Hopf 
algebras $\mathcal{U}_K$ and $\mathcal{A}_K$ 
over the $\mathbb{C}$-algebra $K$ respectively, cf. \S 2.  We denote 
$ab=m_{\mathcal{A}}(a\otimes b)$ for the multiplication 
of the two elements $a,b\in \mathcal{A}_K$ in the $K$-algebra $\mathcal{A}_K$. 
The extended 
comultiplication $\Delta_{\mathcal{A}}$ can be viewed as map
$\Delta_{\mathcal{A}}: \mathcal{A}_K\rightarrow 
\mathcal{A}_K\otimes_K\mathcal{A}_K$ as well as map
$\Delta_{\mathcal{A}}: \mathcal{A}_K\rightarrow
K\otimes \mathcal{A}\otimes \mathcal{A}$ via the canonical identification
$\mathcal{A}_K\otimes_K\mathcal{A}_K\simeq \mathcal{A}_K^{\otimes 2}=
K\otimes \mathcal{A}\otimes \mathcal{A}$.
We extend the $(\mathcal{U}-\mathcal{U})$-bimodule structure 
$K$-linearly to a $(\mathcal{U}_K-\mathcal{U}_K)$-bimodule 
structure on $\mathcal{A}_K$. The actions are given by the formula
\eqref{actiondot}, with the Hopf pairing $\langle \cdot,\cdot\rangle$
extended $K$-bilinearly to a Hopf pairing between $\mathcal{U}_K$
and $\mathcal{A}_K$. Similarly, 
we have a componentwise extension of the bimodule structure on
$\mathcal{A}_K$ to a 
$(\mathcal{U}^{\otimes n}_K-\mathcal{U}^{\otimes n}_K)$-bimodule
structure on $\mathcal{A}^{\otimes n}_K$. 

In the following lemma we list some basic properties of the
regular $\mathcal{U}_K$-action on $\mathcal{A}_K$. The straightforward
proof is left to the reader.

\begin{lem}\label{groundactioncompatible}
The $(\mathcal{U}_{K}-\mathcal{U}_{K})$-bimodule
structure on $\mathcal{A}_K$ has the following properties:
\begin{equation}\label{compatibilityrelations}
\begin{split}
u\cdot 1_{\mathcal{A}}\cdot v&=\epsilon_{\mathcal{U}}(u)
\epsilon_{\mathcal{U}}(v)1_{\mathcal{A}},\\
u\cdot m_{\mathcal{A}}\bigl(a\otimes_Kb\bigr)\cdot v&=
m_{\mathcal{A}}\bigl(\Delta_{\mathcal{U}}(u)\cdot (a\otimes_Kb)
\cdot \Delta_{\mathcal{U}}(v)\bigr),\\
\Delta_{\mathcal{A}}(u\cdot a\cdot v)&=
u_2\cdot \Delta_{\mathcal{A}}(a)\cdot v_1,\\
\Delta_{\mathcal{A}}(a)\cdot u_2&=u_1\cdot \Delta_{\mathcal{A}}(a),\\
\epsilon_{\mathcal{A}}(u\cdot a)&=
\epsilon_{\mathcal{A}}(a\cdot u),\\
u\cdot S_{\mathcal{A}}(a)\cdot v&=S_{\mathcal{A}}
(S_{\mathcal{U}}(v)\cdot a\cdot S_{\mathcal{U}}(u))
\end{split}
\end{equation}
for $u,v\in \mathcal{U}_K$ and $a\in\mathcal{A}_K$.
\end{lem}

We assume from now on that $\mathcal{A}$ is $(T-T)$-semisimple,
\begin{equation}\label{bigradingnormal}
\mathcal{A}=\bigoplus_{\alpha,\beta\in\widehat{T}}
\mathcal{A}[\alpha,\beta],
\end{equation}
with $\mathcal{A}[\alpha,\beta]$ consisting of elements $a\in \mathcal{A}$
satisfying $t\cdot a=t^{\beta}a$ and $a\cdot t=t^{\alpha}a$ 
for all $t\in T$. Note that $1_{\mathcal{A}}
\in \mathcal{A}[0,0]$ by the first equality of \eqref{compatibilityrelations}.
The direct sum
decomposition \eqref{bigradingnormal} defines a 
$\widehat{T}$-bigrading of $\mathcal{A}$ due to the second equality
of \eqref{compatibilityrelations}. 

Note that the primitive idempotents 
$\pi_\alpha\in\mathcal{U}$ ($\alpha\in\widehat{T}$)
of $T$ act on $\mathcal{A}$ by
\[\pi_{\alpha}\cdot a\cdot \pi_{\beta}=
\delta_{\alpha,\delta}\,
\delta_{\beta,\gamma}\,a,\qquad
a\in \mathcal{A}_K[\gamma,\delta]
\]
and 
that the $\widehat{T}$-bigrading of $\mathcal{A}$ is compatible
with the $\widehat{T}$-grading of $\mathcal{U}$,
\[\mathcal{U}[\alpha]\cdot \mathcal{A}[\beta,\gamma]\subseteq
\mathcal{A}[\beta,\alpha+\gamma],\qquad
\mathcal{A}[\beta,\gamma]\cdot \mathcal{U}[\alpha]\subseteq
\mathcal{A}[\beta-\alpha,\gamma].
\]
Lemma \ref{groundactioncompatible} implies 
that the $\widehat{T}$-bigrading of $\mathcal{A}_K$ is compatible with
the Hopf-algebra maps of $\mathcal{A}_K$.
\begin{cor}\label{cocompatible}
Let $\alpha,\beta\in \widehat{T}$. Then
\begin{equation*}
\begin{split}
\Delta_{\mathcal{A}}\bigl(\mathcal{A}_K[\alpha,\beta]\bigr)&\subseteq
\bigoplus_{\gamma\in\widehat{T}}\mathcal{A}_K[\alpha,\gamma]\otimes_K 
\mathcal{A}_K[\gamma,\beta],\\
\epsilon_{\mathcal{A}}\bigl(\mathcal{A}_K[\alpha,\beta]\bigr)&=\{0\}\,
\,\,\hbox{ unless } \, \alpha=\beta,\\
S_{\mathcal{A}}\bigl(\mathcal{A}_K[\alpha,\beta]\bigr)&\subseteq
\mathcal{A}_K[-\beta,-\alpha].
\end{split}
\end{equation*}
\end{cor}
The straightforward proof of Corollary \ref{cocompatible}
is left to the reader.


\subsection{The dynamical quantum group $\mathcal{A}^J$}
The constructions in this subsection are motivated by the 
dynamical quantum group constructions 
of Etingof and Varchenko \cite[\S 4]{EV2} and Etingof and 
Nikshych \cite[\S 4.3]{EN1}. Since the proofs in this subsection
are quite analogous to the ones in \cite[\S 4]{EV2}, 
we only indicate their main steps.

We keep the conventions and notations of the previous subsection.
We fix a dynamical twist $J(\lambda)$ for $\mathcal{U}$ with respect to $T$.
\begin{lem}\label{twisthalg}
The $\mathbb{C}$-vectorspace $\mathcal{A}_K$ is a
$T$-algebra with multiplication
\[m^J(a(\lambda,\mu)\otimes b(\lambda,\mu))=
m_{\mathcal{A}}\bigl(J(\mu)\cdot \bigl(a(\lambda+\alpha,\mu+\beta)\otimes
b(\lambda,\mu)\bigr)\cdot J^{-1}(\lambda)\bigr)
\]
for $a(\lambda,\mu)\in\mathcal{A}_K$ and 
$b(\lambda,\mu)\in\mathcal{A}_K[\alpha,\beta]$, 
with unit $1_{\mathcal{A}}$, with
$(\alpha,\beta)$-bigraded pieces $\mathcal{A}_K[\alpha,\beta]$ 
\textup{(}$\alpha, \beta\in\widehat{T}$\textup{)} and with moment maps 
\[\mu_l(f)=f(\lambda)1_{\mathcal{A}},\qquad
\mu_r(f)=f(\mu)1_{\mathcal{A}}
\]
for $f\in F$, where \textup{(}recall\textup{)}
$f(\lambda)=f\otimes 1_F\in K$ and $f(\mu)=1_F\otimes f\in K$.
\end{lem}
\begin{proof}
The fact that $J(\lambda)$ is of zero $T$-weight implies that
$\mathcal{A}_K=\oplus_{\alpha,\beta}\mathcal{A}_K[\alpha,\beta]$
defines a $\widehat{T}$-bigrading with respect to the new multiplication
$m^J$. The second line of \eqref{dyncond}
implies the associativity of $m^J$. 
The first line of \eqref{dyncond} implies that $1_{\mathcal{A}}$
is the unit element with respect to $m^J$. 
The axioms for the moment maps are straightforward.
\end{proof}
We call $m^J$ the {\it $J$-twisted dynamical multiplication on} 
$\mathcal{A}_K$.
We write $\mathcal{A}^J$ for the $\mathbb{C}$-vectorspace $\mathcal{A}_K$,
viewed as $T$-algebra by Lemma \ref{twisthalg}. 
The following proposition provides the link between the $T$-algebra 
$\mathcal{A}^J$
and the Faddeev-Reshetikhin-Takhtajan (FRT) type construction
of dynamical quantum groups.
\begin{prop}\label{FRT}
Suppose $\mathcal{U}$ is quasi-triangular with universal $R$-matrix
$\mathcal{R}$. Let $\mathcal{R}^J(\lambda)$ \textup{(}see
\eqref{Rdyn}\textup{)} be the corresponding dynamical universal
$R$-matrix. Then
\[m^J\bigl(\mathcal{R}^J(\mu)\cdot (a\otimes b)\bigr)=
m^J\bigl((b\otimes a)\cdot \mathcal{R}_{21}^J(\lambda)\bigr),\qquad
\forall\,a,b\in\mathcal{A}\subset \mathcal{A}^J,
\]
where we use the convention that 
the $\mu$ and $\lambda$ dependence of the action of 
the dynamical universal
$R$-matrices end up in the second tensor component. In other words,
\[m_{\mathcal{A}}\bigl(J(\mu)\mathcal{R}^J(\mu)\cdot (a\otimes b)\cdot
J^{-1}(\lambda)\bigr)=
m_{\mathcal{A}}\bigl(J(\mu)\cdot (b\otimes a)\cdot 
\mathcal{R}_{21}^J(\lambda)J^{-1}(\lambda)\bigr),\qquad
\forall\,a,b\in\mathcal{A}.
\]
\end{prop}
\begin{proof}
This follows directly from the well known FRT type commutation relations
\begin{equation}\label{FRTclassical}
m_{\mathcal{A}}\bigl(\mathcal{R}_{21}\cdot (b\otimes a)\bigr)=
m_{\mathcal{A}}\bigl((a\otimes b)\cdot \mathcal{R}\bigr)
\end{equation}
for $a,b\in\mathcal{A}$.
\end{proof}

Let $\mathcal{A}_K\widehat{\otimes}_F\mathcal{A}_K$ be the 
$\mathbb{C}$-linear vector space defined by taking the tensor
product over $F$ with respect to the balancing condition
$(f(\mu)a)\widehat{\otimes}_Fb=a\widehat{\otimes}_F(f(\lambda)b)$
for $f\in F$ and $a,b\in\mathcal{A}_K$. The restricted tensor product 
$\mathcal{A}^J\widetilde{\otimes}\mathcal{A}^J$
naturally identifies as $\mathbb{C}$-vectorspace with the subspace
\[\bigoplus_{\alpha,\beta,\gamma\in\widehat{T}}
\mathcal{A}_K[\alpha,\gamma]\widehat{\otimes}_F\mathcal{A}_K[\gamma,\beta]
\]
of $\mathcal{A}_K\widehat{\otimes}_F\mathcal{A}_K$. 
Let
\[\pi: K\otimes \mathcal{A}\otimes \mathcal{A}\rightarrow
\mathcal{A}_K\widehat{\otimes}_F\mathcal{A}_K
\]
be the unique $\mathbb{C}$-linear map
satisfying $\pi\bigl(f(\lambda)\otimes g(\mu)\otimes a\otimes b\bigr)=
(f(\lambda)a)\widehat{\otimes}_F (g(\mu)b)$ for $f,g\in F$ and 
$a,b\in \mathcal{A}$. 

\begin{prop}\label{hiodJ}
The $T$-algebra $\mathcal{A}^J$ is a $T$-bialgebroid
with coalgebroid maps $\Delta:\mathcal{A}^J\rightarrow
\mathcal{A}^J\widetilde{\otimes}\mathcal{A}^J$ 
and $\epsilon:\mathcal{A}^J\rightarrow I$ defined by
\begin{equation*}
\begin{split}
\Delta(a)&=\pi\bigl(\Delta_{\mathcal{A}}(a)\bigr),\\
\epsilon(a)&=
T_{-\alpha}\bigl(m_F(\epsilon_{\mathcal{A}}(a))\bigr)T_{-\alpha}
\end{split}
\end{equation*}
for $a\in \mathcal{A}_K[\alpha,\beta]$,
where $m_F: K=F\otimes F\rightarrow F$,
$m_F(f(\lambda)\otimes g(\mu))=f(\lambda)g(\lambda)$ 
is the multiplication map of the 
$\mathbb{C}$-algebra $F$.
\end{prop}
\begin{proof}
Corollary \ref{cocompatible} implies that the image of $\Delta$ is contained
in $\mathcal{A}^J\widetilde{\otimes}\mathcal{A}^J$ and that
$\Delta$ and $\epsilon$ preserve the $\widehat{T}$-bigrading.
The maps $\epsilon$ and $\Delta$ are clearly compatible with the moment maps.
To prove that $\epsilon: \mathcal{A}^J\rightarrow I$ and
$\Delta: \mathcal{A}^J\rightarrow \mathcal{A}^J\widetilde{\otimes}
\mathcal{A}^J$ are morphisms of $T$-algebras
it thus remains to show that they are algebra homomorphisms.
This follows from the compatibility of the 
counit $\epsilon_{\mathcal{A}}$
and comultiplication $\Delta_{\mathcal{A}}$ with the 
$(\mathcal{U}_K-\mathcal{U}_K)$-action
on $\mathcal{A}_K$, see Lemma \ref{groundactioncompatible}.
The remainder of the proof is straightforward.
\end{proof}

In the following proposition we define an antipode $S^J$ for 
the $T$-bialgebroid $\mathcal{A}^J$ using the two zero $T$-weighted
elements $K^J$ and $Q^J$ of $\mathcal{U}_K$ associated to $J$, 
see Lemma \ref{KQ}. 
\begin{prop}
The $T$-bialgebroid $\mathcal{A}^J$ is a $T$-Hopf
algebroid with antipode $S^J: \mathcal{A}^J\rightarrow \mathcal{A}^J$ defined
by
\begin{equation}\label{antipode1}
S^J(a(\lambda,\mu))=S_{\mathcal{A}}\bigl(K^J(\lambda-\beta)\cdot
a(\mu-\alpha,\lambda-\beta)\cdot Q^J(\mu)\bigr),\qquad
a(\lambda,\mu)\in\mathcal{A}_K[\alpha,\beta].
\end{equation}
In particular,
\begin{equation}\label{antipode2}
S^J(a)=S_{\mathcal{A}}\bigl(K^J(\lambda-h)\cdot a\cdot Q^J(\mu)\bigr)
\qquad \forall \, a\in\mathcal{A}.
\end{equation}
\end{prop}
\begin{proof}
Lemma \ref{groundactioncompatible} and the fact that $K^J$ and $Q^J$ are
of zero $T$-weight imply that
\[
S^J\bigl(\mathcal{A}_K[\alpha,\beta]\bigr)\subseteq 
\mathcal{A}_K[-\beta,-\alpha].
\]
A straightforward computation then shows that
the linear map $S^J: \mathcal{A}^J\rightarrow \mathcal{A}^J$
defined by \eqref{antipode1} satisfies the required compatibility conditions
with respect to the moment maps of $\mathcal{A}^J$. 
Hence it suffices to prove the antipode identities for
$a\in \mathcal{A}\subset \mathcal{A}^J$. The required antipode 
identities then reduce to
\[m^J\bigl(\hbox{Id}_{\mathcal{A}^J}\otimes S^J\bigr)\bigl(\Delta(a)\bigr)=
\epsilon_{\mathcal{A}}(a)1_{\mathcal{A}}=
m^J\bigl(S^J\otimes \hbox{Id}_{\mathcal{A}^J}\bigr)
\bigl(\Delta(a)\bigr),
\]
which can be proven by direct computations using 
Lemma \ref{groundactioncompatible}, the antipode 
axioms for $\mathcal{U}$ and \eqref{inverse}.
\end{proof}
\begin{defi}
We call the $T$-Hopf algebroid 
$(\mathcal{A}^J=\oplus_{\alpha,\beta}\mathcal{A}_K[\alpha,\beta],
m^J,1_{\mathcal{A}},\Delta,\epsilon,S^J)$
the dynamical quantum group associated to
the five-tuple $(\mathcal{U},\mathcal{A},\langle \cdot,\cdot\rangle,
T,J(\lambda))$.
\end{defi}
\begin{rema}\label{basic}
Associated to the four-tuple
$(\mathcal{U},\mathcal{A},\langle \cdot,\cdot\rangle,T)$
we always have the trivial dynamical quantum group
$\mathcal{A}^1$, whose associated dynamical twist is the unit element
$1=1_{\mathcal{U}^{\otimes 2}_K}\in\mathcal{U}^{\otimes 2}_K$. 
The multiplication and antipode of $\mathcal{A}^1$ are
\begin{equation*}
\begin{split}
&m^{1}\bigl(a(\lambda,\mu)\otimes
b(\lambda,\mu)\bigr)=m_{\mathcal{A}}\bigl(a(\lambda+\alpha,\mu+\beta)
\otimes b(\lambda,\mu)\bigr),\\
&S^{1}\bigl(b(\lambda,\mu)\bigr)=
S_{\mathcal{A}}\bigl(b(\mu-\alpha,\lambda-\beta)\bigr)
\end{split}
\end{equation*}
for $a(\lambda,\mu)\in \mathcal{A}^{1}$ 
and $b(\lambda,\mu)\in\mathcal{A}_K[\alpha,\beta]$, which are
the trivial dynamical extensions of the multiplication and antipode
of $\mathcal{A}$. 
\end{rema}

\subsection{Gauge equivalent dynamical quantum groups}

Let $x(\lambda)\in\mathcal{U}_{F_1}$ be a gauge
transformation and 
$J(\lambda)\in\mathcal{U}^{\otimes 2}_{F_1}$ a dynamical
twist with respect to $T$. 
Recall that the corresponding gauged dynamical twist is given by
\[J_x(\lambda)=
\Delta_{\mathcal{U}}(x(\lambda))
J(\lambda)x_2^{-1}(\lambda)x_1^{-1}(\lambda+h^{(2)})
\in\mathcal{U}^{\otimes 2}_{F_1}.
\]
\begin{prop}\label{gaugeA}
The dynamical quantum groups $\mathcal{A}^J$ and $\mathcal{A}^{J_x}$
are isomorphic. The corresponding $T$-Hopf algebroid isomorphism 
$\phi_x: \mathcal{A}^J\rightarrow \mathcal{A}^{J_x}$ is
\[\phi_x(a)=x(\mu)\cdot a\cdot x^{-1}(\lambda),\qquad a\in \mathcal{A}^J.
\]
\end{prop}
\begin{proof}
Since $x(\lambda)$ is of zero weight, we have 
$\phi_x(\mathcal{A}_K[\alpha,\beta])=
\mathcal{A}_K[\alpha,\beta]$ for $\alpha,\beta\in\widehat{T}$. The fact that
$\epsilon_{\mathcal{U}}(x(\lambda))=1$ implies
$\phi_x(1_{\mathcal{A}})=1_{\mathcal{A}}$. It follows now directly that
$\phi_x$ respects the moment maps. The map $\phi_x$ is an algebra 
homomorphism since
for $a\in \mathcal{A}$ and $b\in \mathcal{A}[\alpha,\beta]$,
\begin{equation*}
\begin{split}
\phi_x(m^J(a\otimes b))&=
x(\mu)\cdot\bigl(m_{\mathcal{A}}(J(\mu)\cdot(a\otimes b)
\cdot J^{-1}(\lambda))\bigr)\cdot x^{-1}(\lambda)\\
&=m_{\mathcal{A}}\bigl(\Delta_{\mathcal{U}}(x(\mu))J(\mu)\cdot (a\otimes b)
\cdot J^{-1}(\lambda)\Delta_{\mathcal{U}}(x^{-1}(\lambda))\bigr)\\
&=m_{\mathcal{A}}\bigl(J_x(\mu)\cdot\bigl(x(\mu+\beta)\cdot 
a\cdot x^{-1}(\lambda+\alpha)\otimes
x(\mu)\cdot b\cdot x^{-1}(\lambda)\bigr)\cdot J_x^{-1}(\lambda)\bigr)\\
&=m^{J_x}\bigl(\phi_x(a)\otimes\phi_x(b)\bigr),
\end{split}
\end{equation*}
hence $\phi_x$ is 
an isomorphism of $T$-algebras.
The proof that $\phi_x$ is an isomorphism 
of $T$-Hopf algebroids follows
from Lemma \ref{groundactioncompatible} by 
direct computations. We give the computations
for the comultiplication and the antipode.
For the compatibility of $\phi_x$ with respect to the
comultiplication we compute for $a\in\mathcal{A}$,
\begin{equation*}
\begin{split}
\bigl(\phi_x\widetilde{\otimes}\phi_x\bigr)(\Delta(a))&=
\sum \bigl(x(\mu)\cdot a_{(1)}\cdot x^{-1}(\lambda)\bigr)\otimes_F
\bigl(x(\mu)\cdot a_{(2)}\cdot x^{-1}(\lambda)\bigr)\\
&=\sum \bigl(a_{(1)}\cdot x^{-1}(\lambda)\bigr)\otimes_F
\bigl(x(\mu)\cdot a_{(2)}\cdot x(\lambda)x^{-1}(\lambda)\bigr)\\
&=\sum \bigl(a_{(1)}\cdot x^{-1}(\lambda)\bigr)\otimes_F\bigl(x(\mu)
\cdot a_{(2)}\bigr)\\
&=\bigl(\pi\circ
\Delta_{\mathcal{A}}\bigr)(x(\mu)\cdot a\cdot x^{-1}(\lambda))\\
&=\Delta(\phi_x(a)).
\end{split}
\end{equation*}
For the compatibility of $\phi_x$ with respect to the
antipode we first note that
\begin{equation*}
\begin{split}
K^{J_x}(\lambda)&=S_{\mathcal{U}}(x^{-1}(\lambda-h))K^J(\lambda)
x^{-1}(\lambda),\\
Q^{J_x}(\lambda)&=x(\lambda-h)Q^J(\lambda)S_{\mathcal{U}}(x(\lambda)).
\end{split}
\end{equation*}
For $a\in \mathcal{A}[\alpha,\beta]$ we then have
\begin{equation*}
\begin{split}
\phi_x(S^J(a))&=x(\mu)\cdot S_{\mathcal{A}}\bigl(K^J(\lambda-\beta)\cdot 
a\cdot Q^J(\mu)\bigr)\cdot x^{-1}(\lambda)\\
&=S_{\mathcal{A}}\bigl(S_{\mathcal{U}}(x^{-1}(\lambda))
K^J(\lambda-\beta)\cdot a\cdot Q^J(\mu)S_{\mathcal{U}}(x(\mu))\bigr)\\
&=S_{\mathcal{A}}\bigl(K^{J_x}(\lambda-\beta)x(\lambda-\beta)\cdot a
\cdot x^{-1}(\mu-\alpha)Q^{J_x}(\mu)\bigr)\\
&=S^{J_x}\bigl(\phi_x(a)),
\end{split}
\end{equation*}
as desired.
\end{proof}

\subsection{The dynamical quantum group associated to a 
vertex-IRF transformation}\label{twoways}

For a vertex-IRF transformation $x(\lambda)\in \mathcal{U}_{F_1}$ 
with respect to $T$ we define
\begin{equation}\label{xgrading}
\mathcal{A}^x_{\alpha\beta}=x^{-1}(\mu)\cdot 
\mathcal{A}_K[\alpha,\beta]\cdot x(\lambda),
\qquad \alpha,\beta\in\widehat{T}.
\end{equation}
Note that in 
general $\mathcal{A}^x_{\alpha\beta}\not=\mathcal{A}_K[\alpha,\beta]$
since $x(\lambda)$ is not necessarily of zero $T$-weight.
We define a $T$-Hopf algebroid structure 
on $\mathcal{A}_K$ such that $\mathcal{A}_K=
\oplus_{\alpha,\beta}\mathcal{A}^x_{\alpha\beta}$ is the associated
$\widehat{T}$-bigrading.
Recall the dynamical twist
\[j_x(\lambda)=
\Delta_{\mathcal{U}}(x(\lambda))x_2^{-1}(\lambda)x_1^{-1}(\lambda+h^{(2)})
\]
associated to the vertex-IRF transformation $x(\lambda)$.
\begin{thm}\label{twoway}
Let $x(\lambda)\in\mathcal{U}_{F_1}$ be a vertex-IRF 
transformation with respect to $T$.

{\bf a.}
The $\mathbb{C}$-vectorspace $\mathcal{A}_K$ is a $T$-Hopf 
algebroid
with multiplication
\[m^x\bigl(a(\lambda,\mu)\otimes b(\lambda,\mu)\bigr)=
m_{\mathcal{A}}\bigl(a(\lambda+\alpha,\mu+\beta)\otimes 
b(\lambda,\mu)\bigr),\qquad a\in\mathcal{A},\,\,
b(\lambda,\mu)\in\mathcal{A}^x_{\alpha\beta},
\]
with unit $1_{\mathcal{A}}$, 
with $(\alpha,\beta)$-bigraded pieces $\mathcal{A}_{\alpha\beta}^x$
for $\alpha,\beta\in\widehat{T}$, with moment maps
\[\mu_l(f)=f(\lambda)1_{\mathcal{A}},
\qquad \mu_r(f)=f(\mu)1_{\mathcal{A}},\qquad f\in F,
\]
with comultiplication $\Delta=\pi\circ\Delta_{\mathcal{A}}$, with counit
\[\epsilon^x(a)=
T_{-\alpha}\bigl(m_F(\epsilon_{\mathcal{A}}(a))\bigr)T_{-\alpha},
\qquad a\in\mathcal{A}^x_{\alpha\beta}
\]
and with antipode
\[S^x(a(\lambda,\mu))=
S_{\mathcal{A}}\bigl(a(\mu-\alpha,\lambda-\beta)\bigr),
\qquad a(\lambda,\mu)\in \mathcal{A}^x_{\alpha\beta}.
\]
We write $\mathcal{A}^x$ for 
$\mathcal{A}_K$ viewed as $T$-Hopf algebroid
in this way.

{\bf b.} The map $\phi_x: \mathcal{A}^x\rightarrow \mathcal{A}^{j_x}$ 
defined by
\[\phi_x(a)=x(\mu)\cdot a\cdot x^{-1}(\lambda),\qquad a\in\mathcal{A}^x
\]
is an isomorphism of $T$-Hopf algebroids.
\end{thm}
\begin{proof}
There is a unique $T$-Hopf algebroid structure on $\mathcal{A}_K$
turning the $K$-linear isomorphism
$\phi_x: \mathcal{A}_K\rightarrow \mathcal{A}^{j_x}$, defined by
$\phi_x(a)=x(\mu)\cdot a\cdot x^{-1}(\lambda)$, into an isomorphism
of $T$-Hopf algebroids. A direct computation, which is 
similar to the proof of Proposition \ref{gaugeA} for the special case that
$J(\lambda)=1$ is the trivial dynamical twist,
proves that the resulting $T$-Hopf algebroid structure 
on $\mathcal{A}_K$
is defined by the explicit formulas as given in part {\bf a} of the theorem.
\end{proof}
For $x(\lambda)$ a gauge 
transformation with respect to $T$
we have $\mathcal{A}^1=\mathcal{A}^{x}\simeq
\mathcal{A}^{j_x}$ as $T$-Hopf algebroids, with $\mathcal{A}^1$
the trivial dynamical quantum group associated to
$(\mathcal{U},\mathcal{A},\langle \cdot,\cdot\rangle,T)$, 
cf. Remark \ref{basic}. 
In case of a vertex-IRF transformation $x(\lambda)$ the 
bigrading of $\mathcal{A}^x$ is a nontrivial twisting of the
trivial bigrading of $\mathcal{A}^1$,
but the remaining $T$-Hopf algebroid structures of
$\mathcal{A}^x$ have the same untwisted
form as the $T$-Hopf algebroid structures of 
$\mathcal{A}^{1}$.


\section{Askey-Wilson polynomials and the $\hbox{SL}(2)$ 
(dynamical) quantum group}

The work of Babelon \cite{B} (see also \cite{BBB})
implies that the trigonometric $\hbox{SL}(2;\mathbb{C})$ dynamical
quantum group arises from a vertex-IRF transformation.
In this section we describe the two different realizations of the 
$\hbox{SL}(2;\mathbb{C})$ dynamical quantum group 
(due to Theorem \ref{twoway}), 
leading to an intrinsic link between 
the work of Koelink, Rosengren
\cite{KR} on the trigonometric $\hbox{SL}(2;\mathbb{C})$ dynamical
quantum group and the work of Koornwinder
\cite{Koo1}, Noumi, Mimachi \cite{NM} and Koelink \cite{Koe},
\cite{Koe2} on the $\hbox{SL}(2;\mathbb{C})$ quantum group. 
The important additional property of the vertex-IRF 
transformation which is needed here is Rosengren's \cite{R0} observation  
that the vertex-IRF transformation conjugates the $T$-action
to an action by Koornwinder's \cite{Koo1} twisted primitive elements
(the explicit identification of 
Rosengren's \cite{R0} generalized group element 
with Babelon's \cite{B} vertex-IRF transformation is an unpublished 
observation of Rosengren, see also the introduction of \cite{KR}).

We fix a deformation parameter $q\in\mathbb{C}^*$ which is 
not a root of unity. Let $q^{\frac{1}{2}}$ be a fixed choice
of square root of $q$.
\subsection{The $\hbox{SL}(2)$ quantum group}

Let $\mathcal{U}=\mathcal{U}_q(\mathfrak{s}\mathfrak{l}(2))$ 
be the unital associative
$\mathbb{C}$-algebra with generators $k^{\pm 1}$, $e,f$ and relations
\begin{equation*}
\begin{split}
&kk^{-1}=k^{-1}k=1,\\
&ke=qek,\qquad kf=q^{-1}fk,\\
&ef-fe=\frac{k^2-k^{-2}}{q-q^{-1}}.
\end{split}
\end{equation*}
The algebra $\mathcal{U}$ is a Hopf-algebra with comultiplication
\begin{equation*}
\begin{split}
\Delta_{\mathcal{U}}(k^{\pm 1})&=k^{\pm 1}\otimes k^{\pm 1},\\
\Delta_{\mathcal{U}}(e)&=k\otimes e+e\otimes k^{-1}\\
\Delta_{\mathcal{U}}(f)&=k\otimes f+f\otimes k^{-1},
\end{split}
\end{equation*}
with counit
\[\epsilon_{\mathcal{U}}(k^{\pm 1})=1,\qquad \epsilon_{\mathcal{U}}(e)=
\epsilon_{\mathcal{U}}(f)=0
\]
and with antipode
\[S_{\mathcal{U}}(k^{\pm 1})=k^{\mp 1},\qquad
S_{\mathcal{U}}(e)=-q^{-1}e,\qquad S_{\mathcal{U}}(f)=-qf.
\]
We take $T=\{k^m \, | \, m\in\mathbb{Z}\}\subset G(\mathcal{U})$
as abelian subgroup of 
the group-like elements in $\mathcal{U}$.
The role of the character group $\widehat{T}$ is
taken over by the integers $\mathbb{Z}$, viewed as characters
of $T$ by 
\[\bigl(k^m\bigr)^\alpha=q^{m\alpha/2},\qquad m,\alpha\in\mathbb{Z}.
\] 
Note that $\mathcal{U}$ is $\hbox{ad}(T)$-semisimple 
with spectrum contained in $2\mathbb{Z}$.

The type 1 irreducible, finite dimensional $\mathcal{U}$-representations
are parametrized by $\mathbb{Z}_{\geq 0}$. For $m\in \mathbb{Z}_{\geq 0}$
the corresponding spin $\frac{m}{2}$ representation $V_m$ is a 
$m+1$-dimensional
representation with basis 
$v_r^m$ ($r=-m,2-m,\ldots,m-2,m$) and action
\begin{equation*}
\begin{split}
k^{\pm 1}\,v_r^m&=q^{\pm\frac{r}{2}}v_r^m,\\
e\,v_r^m&=
\frac{\sqrt{\bigl(q^{-\frac{1}{2}(m+r+2)}-q^{\frac{1}{2}(m+r+2)}\bigr)
\bigl(q^{-\frac{1}{2}(m-r)}-q^{\frac{1}{2}(m-r)}\bigr)}}
{q^{-1}-q}\,v_{r+2}^m,\\
f\,v_r^m&=\frac{\sqrt{\bigl(q^{-\frac{1}{2}(m+r)}-q^{\frac{1}{2}(m+r)}\bigr)
\bigl(q^{-\frac{1}{2}(m-r+2)}-q^{\frac{1}{2}(m-r+2)}\bigr)}}
{q^{-1}-q}\,v_{r-2}^m
\end{split}
\end{equation*}
where $v_{m+2}^m=v_{-m-2}^m=0$ by convention. Type 1 refers to the fact that
the modules $V_m$ are $T$-semisimple with spectrum contained in $\mathbb{Z}$.
We denote $\bigl(\cdot,\cdot\bigr): V_m\otimes V_m\rightarrow \mathbb{C}$
for the {\it bilinear} pairing such that 
$\bigl(v_r^m,v_s^m\bigr)=\delta_{r,s}$.
Then 
\begin{equation}
\bigl(Xv,w\bigr)=
\bigl(v,X^\ddagger w\bigr),\qquad X\in \mathcal{U},\,\,v,w\in V_m
\end{equation}
with ${}^\ddagger: 
\mathcal{U}\rightarrow \mathcal{U}$ the unital $\mathbb{C}$-linear
antiinvolution determined by
\[\bigl(k^{\pm 1}\bigr)^\ddagger=k^{\pm 1},\qquad
e^\ddagger=f,\qquad f^\ddagger=e.
\]
Allowing suitable completions, $\mathcal{U}$ is a quasi-triangular
Hopf algebra. We denote 
$\mathcal{R}$ by the corresponding
Drinfeld universal $R$-matrix. 
Its action on $V_1\otimes V_1$ is given by
the matrix
\begin{equation}\label{Runtwist}
\mathcal{R}|_{V_1\otimes V_1}=q^{-\frac{1}{2}}
\left(\begin{matrix} 
q& &0 &0 &0\\
0& &1 &0 &0\\
0& &q-q^{-1} &1 &0\\
0& &0 &0 &q
\end{matrix}\right)
\end{equation}
with respect to the ordered basis $\{v_1^1\otimes v_1^1,
v_1^1\otimes v_{-1}^1,v_{-1}^1\otimes v_1^1, v_{-1}^1\otimes v_{-1}^1\}$.

The quantized function algebra $\mathcal{A}=\mathcal{A}_q[\hbox{SL}(2)]$ is
the Hopf-subalgebra of the Hopf-dual $\mathcal{U}^*$ spanned by the 
matrix coefficients of the finite
dimensional type 1 
$\mathcal{U}$-representations. The Peter-Weyl decomposition
of $\mathcal{A}$ is
\begin{equation*}
\begin{split}
\mathcal{A}&=\bigoplus_{m=0}^{\infty}W(m),\\
W(m)&=\hbox{span}_{\mathbb{C}}\{t_{rs}^m(\cdot)
\,\, | \,\, r,s=-m,2-m,\ldots,m-2,m\},
\end{split}
\end{equation*}
where $t_{rs}^m$ is the matrix coefficient $t_{rs}^m(\cdot)=
\bigl(\,\cdot\,v_s^m,v_r^m\bigr)\in
\mathcal{A}[r,s]$.
The Peter-Weyl decomposition 
is the irreducible decomposition of $\mathcal{A}$, viewed as
$\bigl(\mathcal{U}-\mathcal{U}\bigr)$-bimodule.
Clearly $\mathcal{A}$ is $\bigl(T-T\bigr)$-semisimple. 
Its $(\alpha,\beta)$-bigraded 
piece $\mathcal{A}[\alpha,\beta]$ is nonzero 
when $\alpha$ and $\beta$ are integers having the
same parity.

The Hopf-algebra $\mathcal{A}$ is generated as 
unital $\mathbb{C}$-algebra by the
matrix coefficients of the two-dimensional representation $V_1$, which
we denote by
\begin{equation*}
\left(\begin{matrix} 
\alpha &\beta\\
\gamma &\delta
\end{matrix}\right)=
\left(\begin{matrix}
(\,\cdot\,v_1^1,v_1^1)& &(\,\cdot\,v_{-1}^1,v_1^1)\\
(\,\cdot\,v_1^1,v_{-1}^1)& &(\,\cdot\,v_{-1}^1,v_{-1}^1)
\end{matrix}\right).
\end{equation*}
Note that $\alpha\in\mathcal{A}[1,1]$, $\beta\in \mathcal{A}[1,-1]$,
$\gamma\in\mathcal{A}[-1,1]$ and $\delta\in\mathcal{A}[-1,-1]$. 
The characterizing commutation relations are governed by the
FRT relations \eqref{FRTclassical} for $a,b\in\{\alpha,\beta,\gamma,\delta\}$
and by the determinant relation
\begin{equation}\label{determinant}
\delta\alpha-q^{-1}\beta\gamma=1_{\mathcal{A}}.
\end{equation}
Explicitly, the FRT commutation relations give the relations
\begin{equation*}
\begin{split}
&\alpha\beta=q\beta\alpha,\quad
\alpha\gamma=q\gamma\alpha,\quad
\beta\delta=q\delta\beta,\\
&\gamma\delta=q\delta\gamma,\quad
\beta\gamma=\gamma\beta,\\
&\alpha\delta-\delta\alpha=(q-q^{-1})\beta\gamma.
\end{split}
\end{equation*}
 
Since for $m\in\mathbb{Z}_{\geq 0}$ and $r,s\in\{-m,2-m,\ldots,m-2,m\}$, 
\[\hbox{span}_{\mathbb{C}}\{t_{rs}^m\}=W(m)\cap \mathcal{A}[r,s],
\]
the study of the matrix coefficients $t_{rs}^m$ relates to 
harmonic analysis on the $\hbox{SL}(2;\mathbb{C})$ 
quantum group $\mathcal{A}$ with respect to $T$.
Considering the $t_{rs}^m$ ($r,s=-m,2-m,\ldots,m-2,m$)
as the matrix coefficients of a finite dimensional
$\mathcal{A}$-corepresentation, the study of the
$t_{rs}^m$ relates to harmonic analysis
on $\mathcal{A}$ with respect to the standard quantum analogue 
of the Cartan subalgebra of $\mathfrak{s}\mathfrak{l}(2;\mathbb{C})$. 

A first example relating basic hypergeometric series
to harmonic analysis on quantum groups is the 
expression of the coefficients
$t_{rs}^m$ ($r,s\in\{-m,2-m,\ldots,m-2,m\}$) in terms of little 
$q$-Jacobi polynomials, see \cite{VS}, \cite{KooIndag} and
\cite{M}. As a special case we recall the formula for the 
matrix coefficients $t_{rs}^m$ with integers $r,s,m$ having the same parity
and satisfying $-m\leq r\leq s\leq -r\leq m$, given by
\[t_{rs}^{m}=C_{rs}^m\,\delta^{\frac{-r-s}{2}}\gamma^{\frac{s-r}{2}}\,
{}_2\phi_1\left(\begin{matrix} q^{-r-m},
q^{-r+m-2}\\ q^{s-r+2}\end{matrix};\,q^2,-q\beta\gamma\right)
\]
for some explicit nonzero constant $C_{rs}^m$. Here
\[{}_{r+1}\phi_r\left(\begin{matrix} a_1,\ldots,a_{r+1}\\
b_1,\ldots,b_r\end{matrix};q,z\right)=
\sum_{m=0}^{\infty}\frac{\bigl(a_1;q\bigr)_m\cdots \bigl(a_{r+1};q\bigr)_m}
{\bigl(q;q\bigr)_m\bigl(b_1;q\bigr)_m\cdots \bigl(b_r;q\bigr)_m}\,
z^m,
\]
with 
$\bigl(a;q\bigr)_m=\prod_{j=0}^{m-1}(1-aq^j)$ 
($m\in \mathbb{Z}_{\geq 0}\cup\{\infty\}$) the $q$-shifted factorial, 
is the ${}_{r+1}\phi_r$ basic hypergeometric series, see \cite{GR}.
The element $\beta\gamma$ is an algebraic generator of
the unital $\mathbb{C}$-subalgebra 
$\mathcal{A}[0,0]$ of $\mathcal{A}$, and it is ``quasi-central'' in
$\mathcal{A}$ (it quasi-commutes with the four 
generators $\alpha,\beta,\gamma$ and
$\delta$). The quasi-centrality can be best 
reformulated in dynamical terms:
$q^{\lambda+\mu}\beta\gamma$ is 
a central element of the trivial dynamical quantum group
$\mathcal{A}^{1}$, cf. Remark \ref{basic}.


\subsection{The $\hbox{SL}(2)$ dynamical quantum group}
The four-tuple $(\mathcal{U},\mathcal{A},\langle \cdot,\cdot\rangle,T)$
as constructed in \S 4.1 does not quite fit into the formal algebraic
setup of \S 3 since the abelian group 
$T$ is not finite. The results and constructions
of \S 3 though still hold true in the present setup by interpreting
the action of the idempotents $\pi_\alpha$ ($\alpha\in\mathbb{Z}$)
on $\mathcal{A}_K$ as the projection operators
\[\pi_\alpha\cdot a\cdot \pi_\beta=
\delta_{\alpha,\delta}\delta_{\beta,\gamma}\,a,\qquad
a\in \mathcal{A}[\gamma,\delta].
\]
Furthermore it is convenient to replace the role of the 
function algebra $F$ in \S 3 by the field $F$ of meromorphic 
functions on $\mathbb{C}$ and accordingly we take $K=F\otimes F$.

Besides this formal extension of the setup of \S 3, 
we also need to work with a suitable completion
of the algebra $\mathcal{U}_K$, which we do not specify here in detail.
All explicit formulas given below will have an obvious, functional calculus
type meaning when acting on $\mathcal{A}_K$ via the left or right regular 
action on $\mathcal{A}$. In particular, all infinite sums below become
finite when acting on $\mathcal{A}$ since both $e\in \mathcal{U}$
and $f\in\mathcal{U}$ act locally nilpotently on $\mathcal{A}_K$.

The upshot is that all universal expressions and universal identities in
(a suitable completions of) $\mathcal{U}_K$ given below should be interpreted 
within $\hbox{End}_K(\mathcal{A}_K)$ through the representation maps
of the left and right regular $\mathcal{U}_K$-action on $\mathcal{A}_K$,
and as such the results of \S 3 hold true.

Babelon \cite{B}, see also 
\cite[\S 2]{BBB} and \cite[\S 7]{KT}, considered the element
\begin{equation}\label{x}
x(\lambda)=\sum_{l,m=0}^{\infty}
\frac{\bigl(-q^{-\lambda}\bigr)^{l+m}}
{\bigl(q^{-2\lambda};q^2\bigr)_{l}}
\frac{q^{l^2-l+2lm-2m}(1-q^2)^{l+m}}{\bigl(q^2;q^2\bigr)_l
\bigl(q^2;q^2\bigr)_m}
\bigl(fk^{-1}\bigr)^l\bigl(ek^{-1}\bigr)^m.
\end{equation} 
The element $x(\lambda)$ is directly related to Rosengren's 
\cite[Prop. 3.3]{R0} group element $U_{\lambda\mu}$
by 
\[x(\lambda)=\bigl(U_{q^{-\lambda-1},q^{-\lambda-1}}\bigr)^\dagger,
\]
where we identify $q^{\frac{1}{2}},K^{\pm},X_+,X_-$ in \cite{R0} 
with $q,k^{\pm 1},e,f$ respectively, 
and where $\dagger$ is the $K$-linear 
antiinvolution of $\mathcal{U}_K$ defined by
\[\bigl(k^{\pm 1}\bigr)^\dagger=k^{\pm 1},\qquad
e^\dagger=-q^{\frac{1}{2}}f,\quad f^\dagger=-q^{-\frac{1}{2}}e.
\]
Thus \cite[Prop. 3.5]{R0} implies that $x(\lambda)$ is invertible
(as element of $\hbox{End}_K(\mathcal{A}_K)$) with inverse given by
\begin{equation}\label{xinverse}
\begin{split}
x^{-1}(\lambda)&=\sum_{m=0}^{\infty}\frac{q^{-m\lambda}q^{-m}}
{\bigl(q^2;q^2\bigr)_m}(1-q^2)^m(k^{-1}f)^m\\
&\times \sum_{l=0}^{\infty}
\frac{q^{-l\lambda}q^{l^2-2l}}{\bigl(q^2;q^2\bigr)_l}
(1-q^2)^l(k^{-1}e)^l\frac{\bigl(q^{-2\lambda};q^2\bigr)_{\infty}}
{\bigl(q^{-2(l+\lambda+1)}k^{-4};q^2\bigr)_l
\bigl(q^{-2\lambda}k^{-4};q^2\bigr)_{\infty}}.
\end{split}
\end{equation}
Babelon \cite{B}, see also \cite{BBB}, observed that $x(\lambda)$ is a 
vertex-IRF transformation with respect to $T$. 
The corresponding dynamical twist 
$j_x(\lambda)$, cf. \cite[\S 1]{BBB}, is directly related to the 
universal fusion matrix for $\mathcal{U}$, see \cite{EV3} and \cite{KT}. 
The explicit expressions for $j_x(\lambda)$ and $j_x^{-1}(\lambda)$,
regarded as elements of $\hbox{End}_K(\mathcal{A}_K)^{\otimes 2}$ 
through the left and right regular representation of $\mathcal{U}_K$, are
\begin{equation}\label{j}
\begin{split}
j_x(\lambda)&=\sum_{l=0}^{\infty}(1-q^2)^{2l}\,
\frac{(-1)^lq^{-2l\lambda+2l^2-4l}}{\bigl(q^2;q^2\bigr)_l}
\left(k^{-l}f^l\otimes \frac{1}{\bigl(q^{2(l-1-\lambda)}k^{-4};q^2\bigr)_l}
k^{-3l}e^l\right),\\
j_x^{-1}(\lambda)&=\sum_{l=0}^{\infty}(1-q^2)^{2l}\,
\frac{q^{-2l\lambda+l^2-3l}}{\bigl(q^2;q^2\bigr)_l}
\left(k^{-l}f^l\otimes \frac{1}{\bigl(q^{-2\lambda}k^{-4};q^2\bigr)_l}
k^{-3l}e^l\right),
\end{split}
\end{equation}
see \cite[\S 2]{BBB} and \cite[\S 7]{KT}. 

The correponding universal dynamical 
$R$-matrix $\mathcal{R}^{j_x}(\lambda)$, acting on the
representation space $K\otimes V_1\otimes V_1$, 
can be computed explicitly
using \eqref{Rdynx} and \eqref{Runtwist}.  
With respect to
the ordered basis 
$\{v_1^1\otimes v_1^1,v_1^1\otimes v_{-1}^1, v_{-1}^1\otimes v_1^1,
v_{-1}^1\otimes v_{-1}^1\}$ it is given by
\begin{equation}\label{Rdynsl2}
\mathcal{R}^{j_x}(\lambda)|_{V_1\otimes V_1}=q^{-\frac{1}{2}}
\left(\begin{matrix}
q& &0 &0 &0\\
0& &1 &\frac{q^{-1}-q}{q^{-2(\lambda+1)}-1} &0\\
0& &\frac{q^{-1}-q}{q^{2(\lambda+1)}-1} 
&\frac{(q^{-2(\lambda+1)}-q^2)(q^{-2(\lambda+1)}-q^{-2})}
{(q^{-2(\lambda+1)}-1)^2} &0\\
0& &0 &0 &q
\end{matrix}\right).
\end{equation}
The dynamical quantum group $\mathcal{A}^{j_x}$ can now directly be related
to Koelink's and Rosengren's 
\cite[Def. 2.4]{KR} trigonometric $\hbox{SL}(2;\mathbb{C})$
dynamical quantum group as follows.
Observe that the dynamical quantum group $\mathcal{A}^{j_x}$ 
is generated as unital algebra
by $\mu_l(F)$, $\mu_r(F)$ and the matrix 
coefficients $\alpha,\beta,\gamma,\delta$
of the two-dimensional $\mathcal{U}$-representation $V_1$. 
By abuse of notation we denote $f(\lambda)$ and $g(\mu)$ for the
elements $\mu_l(f)$ and $\mu_r(g)$ in $\mathcal{A}^{j_x}$. 
The defining commutation
relations (with multiplication denoted by 
$\circ$ to distinguish it from the
multiplication in $\mathcal{A}_K$), are
$f(\lambda)\circ g(\mu)=g(\mu)\circ f(\lambda)$, 
\begin{equation*}
\begin{split}
f(\lambda)\circ\alpha&=\alpha\circ f(\lambda+1),\qquad 
f(\lambda)\circ\beta=\beta\circ f(\lambda+1),\\
f(\lambda)\circ\gamma&=\gamma\circ f(\lambda-1),\qquad 
f(\lambda)\circ\delta=\delta\circ f(\lambda-1),\\
f(\mu)\circ\alpha&=\alpha\circ f(\mu+1),\qquad 
f(\mu)\circ\beta=\beta\circ f(\mu-1),\\
f(\mu)\circ\gamma&=\gamma\circ f(\mu+1),\qquad 
f(\mu)\circ\delta=\delta\circ f(\mu-1),
\end{split}
\end{equation*}
the dynamical FRT commutation relations for 
$\alpha,\beta,\gamma$ and $\delta$
(see Proposition \ref{FRT}), 
and finally the dynamical determinant identity
\[\delta\circ\alpha-q^{-1}\left(\frac{q^2-q^{2(\lambda+1)}}
{1-q^{2(\lambda+1)}}\right)\circ\beta\circ\gamma=
1_{\mathcal{A}},
\]
which is simply the 
determinant identity \eqref{determinant} rewritten in terms of the 
$j_x$-twisted dynamical
multiplication $m^{j_x}$. The dynamical FRT commutation relations can 
be expressed in the
familiar form
\begin{equation*}
\sum_{y,y^\prime}L_{y^\prime \xi^\prime}\circ L_{y\xi}\circ
R_{\eta^\prime \eta}^{y^\prime y}(\lambda)=
\sum_{y,y^\prime}L_{\eta y}\circ L_{\eta^\prime y^\prime}\circ 
R_{y^\prime y}^{\xi^\prime \xi}(\mu)
\end{equation*}
with the indices from $\{\pm 1\}$, where the $L_{\xi\eta}$ are given by
\[
L_{1,1}=\alpha,\quad L_{1,-1}=\beta,\quad L_{-1,1}=\gamma,\quad
L_{-1,-1}=\delta, 
\]
and with the coefficients $R^{\xi\xi^\prime}_{\eta\eta^\prime}(\lambda)$
defined by 
\[\mathcal{R}^{j_x}(\lambda)\bigl(v_\xi^1\otimes v_{\xi^\prime}^1\bigr)=
\sum_{\eta,\eta^\prime} 
R^{-\xi,-\xi^\prime}_{-\eta,-\eta^\prime}(\lambda)v_\eta^1\otimes
v_{\eta^\prime}^1.
\]
The equivalence
with the $\hbox{SL}(2;\mathbb{C})$ dynamical quantum 
group of Koelink and Rosengren \cite{KR} now follows
by identifying the
generators $(f(\lambda),g(\mu),\alpha,\beta,\gamma,\delta)$ 
in \cite[Def. 2.4]{KR}
with $(f(-\lambda-2),g(-\mu-2),\delta,\gamma,\beta,\alpha)$.
 
Harmonic analysis on the dynamical quantum group $\mathcal{A}^{j_x}$
with respect to the standard quantum analogue of the Cartan subalgebra
of $\mathfrak{s}\mathfrak{l}(2;\mathbb{C})$
still amounts to the study of the matrix 
coefficients $t_{rs}^m$, now viewed
as matrix coefficients of a tempered corepresentation of 
the dynamical quantum group $\mathcal{A}^{j_x}$ (see \cite[\S 3]{KR}).
In analogy with the harmonic analysis of the ordinary quantum group
$\mathcal{A}$, we now express the matrix coefficient 
$t_{rs}^m\in W(m)\cap \mathcal{A}_K[r,s]$
in terms of a central element $\Xi\in \mathcal{A}^{j_x}$ which, together
with $\mu_l(F)$ and $\mu_r(F)$, generate $\mathcal{A}_K[0,0]$ as unital
$\mathbb{C}$-algebra. The element $\Xi$ is given explicitly by
\begin{equation}
\Xi=q^{-\lambda+\mu+1}+q^{\lambda-\mu-1}-
q^{\lambda+\mu+2}\bigl(1-q^{-2\lambda}\bigr)
\bigl(1-q^{-2(\mu+2)}\bigr)\circ\beta\circ\gamma,
\end{equation}
see \cite[Lem. 3.3]{KR}. 
The matrix coefficient $t_{rs}^m$
for integers $r,s,m$ having the same parity
and satisfying $-m\leq r\leq s\leq -r\leq m$ is then given by
\begin{equation}\label{AWcase}
t_{rs}^{m}=C_{rs}^m(\lambda,\mu)\circ
\delta^{\circ\frac{-r-s}{2}}\circ\gamma^{\circ\frac{s-r}{2}}\circ\,
p_{\frac{m+r}{2}}\bigl(\Xi;\, q^{\lambda-\mu+1}, q^{-\lambda+\mu+1-r+s},
q^{-\lambda-\mu-1-r-s}, q^{\lambda+\mu+3};q^2\bigr)
\end{equation}
for some explicit nonzero meromorphic function $C_{rs}^m(\lambda,\mu)$,
with $p_n$ the Askey-Wilson polynomial of degree $n$,
\[p_n(z+z^{-1};a,b,c,d;q)={}_4\phi_3
\left(\begin{matrix} q^{-n},
abcdq^{n-1}, az, az^{-1}\\ 
ab, ac, ad\end{matrix};\,q,q\right),
\]
see \cite[Thm. 3.5]{KR}. Here we used the notation $\delta^{\circ m}=
\delta\circ\cdots\circ \delta$ ($m$ times) and similarly for $\gamma$. 
The polynomial expression in $\Xi$ in formula \eqref{AWcase} has the 
obvious interpretation
as element in the commutative subalgebra $\mathcal{A}_K[0,0]$ of the 
dynamical quantum group $\mathcal{A}^{j_x}$.

\subsection{The alternative realization $\mathcal{A}^x$}
We define a Cartan type element $X(\lambda)\in \mathcal{U}_{F_{1}}$
and a twisted primitive element $Y(\lambda)\in \mathcal{U}_{F_1}$
by  
\begin{equation*}
\begin{split}
X(\lambda)&=\frac{q^{-\lambda-1}(k^{-2}-1)+q^{\lambda+1}(k^2-1)}
{q-q^{-1}},\\
Y(\lambda)&=fk-ek+\left(\frac{q^{-\lambda-1}+q^{\lambda+1}}
{q-q^{-1}}\right)(k^2-1).
\end{split}
\end{equation*}
The twisted primitive element $Y(\lambda)$ was introduced by Koornwinder
\cite{Koo1} as an one-pa\-ra\-me\-ter family of quantum analogues of 
Lie-algebra generators for nonstandard Cartan subalgebras of 
$\mathfrak{s}\mathfrak{l}(2;\mathbb{C})$.
 
A key property of the vertex-IRF transformation 
$x(\lambda)$ (see \eqref{x}),
proven by Rosengren \cite{R0}, is the fact that
\begin{equation}\label{Rotwist}
x(\lambda)Y(\lambda)\,x^{-1}(\lambda)=X(\lambda)
\end{equation}
viewed as identity in $\hbox{End}_K(\mathcal{A}_K)$ via the left
or right regular representation of $\mathcal{U}_K$, see \cite[(4.8)]{R0}.

We now study the dynamical quantum group $\mathcal{A}^x$ 
using the $T$-Hopf algebroid isomorphism
$\phi_x: \mathcal{A}^x\rightarrow \mathcal{A}^{j_x}$, 
\[\phi_x(a)=x(\mu)\cdot a\cdot x^{-1}(\lambda),\qquad a\in \mathcal{A}^x,
\]
see Theorem \ref{twoway}. As observed in \S \ref{twoways}, the bigrading
of $\mathcal{A}^x$ is nontrivial, but the other $T$-Hopf algebroid
structures have the same untwisted form as the $T$-Hopf algebroid
structures of the trivial dynamical quantum group $\mathcal{A}^1$.
Thus an explicit description of the bigraded pieces 
$\mathcal{A}^x_{\alpha\beta}$ completely clarifies the $T$-Hopf algebroid
structure of $\mathcal{A}^x$.
\begin{prop}
Denote
\[\nu_\alpha(\lambda)=\frac{q^{-\lambda-1}(q^{-\alpha}-1)+
q^{\lambda+1}(q^\alpha-1)}{q-q^{-1}}\in F
\]
for $\alpha\in \mathbb{Z}$. Then
\[\mathcal{A}_{\alpha\beta}^x=\{ a\in \mathcal{A}_K \, | \,
Y(\mu)\cdot a=\nu_\beta(\mu)a,\quad
a\cdot Y(\lambda)=\nu_\alpha(\lambda)a \}
\]
for $\alpha,\beta\in\mathbb{Z}$.
\end{prop}
\begin{proof}
This follows from \eqref{Rotwist} and the observe that
\[\mathcal{A}_K[\alpha,\beta]=\{a\in\mathcal{A}_K \, | \, 
X(\mu)\cdot a=\nu_\beta(\mu)a,\quad a\cdot X(\lambda)=\nu_\alpha(\lambda)a \}.
\]
\end{proof}
Denote for $m\in\mathbb{Z}_{\geq 0}$ and $r,s\in \{-m,2-m,\ldots,
m-2,m\}$,
\[a^m_{rs}(\lambda,\mu)=\phi_x^{-1}(t_{rs}^m)=
x^{-1}(\mu)\cdot t^m_{rs}\cdot x(\lambda),
\]
then $W_K(m)\cap \mathcal{A}_{rs}^x$ is an one-dimensional $K$-subspace
of $\mathcal{A}^x$, spanned by $a^m_{rs}(\lambda,\mu)$. 
Harmonic analysis on $\mathcal{A}$
with respect to Koornwinder's twisted primitive elements precisely amounts
to the study of the matrix coefficients $a^m_{rs}(\lambda,\mu)$. 
Thus the isomorphism $\phi_x$ yields the equivalence between harmonic
analysis on $\mathcal{A}^{j_x}$ with respect to the standard quantum 
Cartan subalgebra in \cite{KR}, and
the harmonic analysis on $\mathcal{A}$ with respect to Koornwinder's
twisted primitive elements as studied in \cite{Koo1}, \cite{NM},
\cite{Koe}. 

To be more concrete, we end this article
by translating the results of the previous subsection
to the twisted primitive picture and linking it to the results
in \cite{Koo1}, \cite{NM}, \cite{Koe} and \cite{Koe2}.
The generators
\begin{equation*}
\begin{split}
\alpha(\lambda,\mu)&=
a^1_{1,1}(\lambda,\mu)=x^{-1}(\mu)\cdot \alpha\cdot x(\lambda),\\
\beta(\lambda,\mu)&=a^1_{1,-1}(\lambda,\mu)=
x^{-1}(\mu)\cdot \beta\cdot x(\lambda),\\
\gamma(\lambda,\mu)&=
a^1_{-1,1}(\lambda,\mu)=x^{-1}(\mu)\cdot \gamma\cdot x(\lambda),\\
\delta(\lambda,\mu)&=
a^1_{-1,-1}(\lambda,\mu)=x^{-1}(\mu)\cdot\delta\cdot x(\lambda)
\end{split}
\end{equation*}
of $\mathcal{A}^x$ can be rewritten in terms of the standard generators
$\alpha,\beta,\gamma$ and $\delta$ of $\mathcal{A}_K$ by
\begin{equation}\label{abcd}
\begin{split}
\alpha(\lambda,\mu)&=\frac{1}{\bigl(1-q^{-2(\mu+1)}\bigr)}
\bigl(\alpha+q^{-\mu-\frac{1}{2}}\beta-q^{-\lambda-\frac{3}{2}}\gamma-
q^{-\lambda-\mu-2}\delta\bigr),\\
\beta(\lambda,\mu)&=\frac{\bigl(1-q^{-2\mu}\bigr)}{\bigl(1-q^{-2(\mu+1)}\bigr)}
\bigl(q^{-\mu-\frac{3}{2}}\alpha+\beta-q^{-\lambda-\mu-3}\gamma
-q^{-\lambda-\frac{3}{2}}\delta\bigr),\\
\gamma(\lambda,\mu)&=\frac{1}{\bigl(1-q^{-2(\mu+1)}\bigr)
\bigl(1-q^{-2\lambda}\bigr)}
\bigl(-q^{-\lambda-\frac{1}{2}}\alpha-q^{-\lambda-\mu-1}\beta
+\gamma+q^{-\mu-\frac{1}{2}}\delta\bigr),\\
\delta(\lambda,\mu)&=\frac{\bigl(1-q^{-2\mu}\bigr)}
{\bigl(1-q^{-2\lambda}\bigr)\bigl(1-q^{-2(\mu+1)}\bigr)}
\bigl(-q^{-\lambda-\mu-2}\alpha-q^{-\lambda-\frac{1}{2}}\beta+
q^{-\mu-\frac{3}{2}}\gamma+\delta\bigr).
\end{split}
\end{equation}
These expressions are easily derived using 
the explicit expression \eqref{x} for the 
vertex-IRF transformation $x(\lambda)$.
Identifying the elements $\alpha,\beta,\gamma,\delta,
A,B,C,D,q^\sigma,q^\tau$ in \cite{Koe} with 
\[
\alpha,q^{\frac{1}{2}}\beta,q^{-\frac{1}{2}}\gamma,\delta,
k,q^{-\frac{1}{2}}e, q^{\frac{1}{2}}f,k^{-1},\sqrt{-1}q^{-\mu-1},
\sqrt{-1}q^{-\lambda-1},
\]
the elements \eqref{abcd} correspond, up to 
$K$-normalization,
to the elements $k^{-1}\cdot\alpha_{\tau,\sigma}$, 
$k^{-1}\cdot\beta_{\tau,\sigma}$, 
$k^{-1}\cdot\gamma_{\tau,\sigma}$
and $k^{-1}\cdot\delta_{\tau,\sigma}$ of 
\cite[Prop. 6.5]{Koe}.

We denote $\rho(\lambda,\mu)=\phi_x^{-1}(\Xi)\in\mathcal{A}^x$, so
\begin{equation}\label{rho}
\rho(\lambda,\mu)=
q^{-\lambda+\mu+1}+q^{\lambda-\mu-1}-q^{\lambda+\mu+2}
\bigl(1-q^{-2\lambda}\bigr)\bigl(1-q^{-2(\mu+2)}\bigr)
\beta(\lambda-1,\mu+1)\gamma(\lambda,\mu).
\end{equation}
By a direct computation using \eqref{abcd}, the element 
$\rho(\lambda,\mu)$
can be explicitly represented as a quadratic 
expression in $\alpha,\beta,\gamma$
and $\delta$. Identifying the generators of \cite{Koe} with ours as
indicated above, $\rho(\lambda,\mu)$
equals $2k^{-1}\cdot \rho_{\tau,\sigma}$, with 
$\rho_{\tau,\sigma}$ the element as defined in \cite[Thm. 5.1]{Koe}
(it was initially written down explicitly in \cite{Koo1}).
The pre-image under $\phi_x$ of the expression \eqref{AWcase} for integers
$r,s,m$ having the same parity and satisfying $-m\leq r\leq s\leq -r\leq m$
yields
\begin{equation*}
\begin{split}
a_{rs}^m(\lambda,\mu)=
D_{rs}^m(\lambda,\mu)&\prod_{i=0}^{-\frac{r}{2}-\frac{s}{2}-1}
\delta(\lambda+r+1+i,\mu+s+1+i)\\
\times&\prod_{j=0}^{-\frac{r}{2}+\frac{s}{2}-1}
\gamma\bigl(\lambda+\frac{r}{2}-\frac{s}{2}+1+j,
\mu-\frac{r}{2}+\frac{s}{2}-1-j\bigr)\\
\times&p_{\frac{m+r}{2}}\bigl(\rho(\lambda,\mu); q^{\lambda-\mu+1},
q^{-\lambda+\mu+1-r+s}, q^{-\lambda-\mu-1-r-s}, q^{\lambda+\mu+3};q^2\bigr)
\end{split}
\end{equation*}
for some non-zero meromorphic function $D_{rs}^m(\lambda,\mu)$,
where $\prod_{i=0}^jb_i$  with $b_i\in\mathcal{A}_K$
equals $1_{\mathcal{A}}$ when $j<0$ and equals $b_0b_1\cdots b_j$ 
when $j\geq 0$.
This formula is in accordance with the 
expressions derived by Koornwinder \cite{Koo1}
(in case $r=s=0$ and $m$ even), and by Noumi, 
Mimachi \cite{NM} and Koelink \cite{Koe}
for arbitrary $r,s,m$, 
compare for instance with 
the expressions in \cite[Cor. 7.8(i)]{Koe}
and \cite[Cor. 6.8]{Koe}.
 
The above dynamical quantum group interpretation of the
harmonic analysis with respect to twisted 
primitive elements  
leads to natural interpretations and 
new proofs of several other known facts. 
For instance, the statement \cite[Prop. 6.5]{Koe} amounts to a 
reformulation of the fact that $\mathcal{A}^x=\oplus_{\alpha,\beta}
\mathcal{A}^x_{\alpha\beta}$ defines a bigrading with respect to the
untwisted dynamical multiplication $m^x$. The factorized form \eqref{rho}
of $\rho(\lambda,\mu)$ is precisely \cite[Prop. 4.1.7]{Koe2}.
The centrality of 
$\rho(\lambda,\mu)$ in the dynamical quantum group $\mathcal{A}^x$
is \cite[Cor. 4.1.8]{Koe2}. The interpretation of the $t_{rs}^m$
as a matrix corepresentation 
of the dynamical quantum group $\mathcal{A}^{j_x}$, 
\begin{equation*}
\begin{split}
\Delta(t_{rs}^m)&=\sum_lt_{rl}^m\otimes_F t_{ls}^m,\\
\epsilon(t_{rs}^m)&=\delta_{r,s}T_{-r},
\end{split}
\end{equation*}
implies via the isomorphism $\phi_x: \mathcal{A}^{x}\rightarrow
\mathcal{A}^{j_x}$ that the $a_{rs}^m(\lambda,\mu)$ define a 
matrix corepresentation of the dynamical quantum group $\mathcal{A}^x$,
\begin{equation*}
\begin{split}
\Delta(a_{rs}^m(\lambda,\mu))&=\sum_la_{rl}^m(\lambda,\mu)\otimes_F
a_{ls}^m(\lambda,\mu),\\
\epsilon^x(a_{rs}^m(\lambda,\mu))&=\delta_{r,s}T_{-r}.
\end{split}
\end{equation*}
The latter formulas directly relate to \cite[Prop. 6.1.1]{Koe2}.



\begin{thebibliography}{XX}
\bibitem[1]{AW} R. Askey, J.A. Wilson, 
{\it Some basic hypergeometric polynomials
that generalize Jacobi polynomials}, Mem. Amer. Math. Soc. {\bf 54}
(1985), no. 319.
\bibitem[2]{B} O. Babelon, {\it Universal exchange algebra for Bloch
waves and Liouville theory}, Comm. Math. Phys.
{\bf 139} (1991), no. 3, pp. 619--643.
\bibitem[3]{BBB} O. Babelon, D. Bernard, E. Billey,
{\it A quasi-Hopf algebra interpretation of quantum $3-j$ and
$6-j$ symbols and difference equations},
Phys. Lett. B {\bf 375} (1996), no. 1-4, pp. 89--97.
\bibitem[4]{EN1} P. Etingof, D. Nikshych, {\it Dynamical quantum groups
at root of 1}, Duke Math. J. {\bf 108} (2001), no. 1, pp. 135--168.
\bibitem[5]{EN2} P. Etingof, D. Nikshych, {\it Vertex-IRF transformations 
and quantization of dynamical r-matrices}, Math. Res. Lett. {\bf 8}
(2001), pp. 331--345.
\bibitem[6]{EV1} P. Etingof, A. Varchenko, {\it Solutions of the quantum
dynamical Yang-Baxter equation and dynamical quantum groups},
Comm. Math. Phys. {\bf 196} (1998), no. 3, pp. 591--640.
\bibitem[7]{EV2} P. Etingof, A. Varchenko, {\it Exchange dynamical quantum
groups}, Comm. Math. Phys. {\bf 205} (1999), no. 1, pp. 19--52. 
\bibitem[8]{EV3} P. Etingof, A. Varchenko, {\it Traces of intertwiners
for quantum groups and difference equations, I}, Duke Math. J. {\bf 104},
(2000), no. 3, pp. 391--432.
\bibitem[9]{F} G. Felder, {\it Conformal field theory and integrable
systems associated to elliptic curves}, Proc. Int. Congress of Math.,
Vol. 1,2 (Z{\"u}rich, 1994), pp. 1247--1255. 
\bibitem[10]{GR} G. Gasper, M. Rahman, {\it Basic Hypergeometric Series},
Encyclopedia of Mathematics and its Applications {\bf 35},
Cambridge Univ. Press (1990).
\bibitem[11]{GN} J.-L. Gervais, A. Neveu, {\it Novel triangle relation
and absence of tachyons in Liouville string field theory},
Nucl. Phys. B {\bf 238} (1984), pp. 125--141.
\bibitem[12]{Koe} E. Koelink, {\it Askey-Wilson polynomials and the quantum
$\hbox{SU}(2)$ group: survey and applications}, Acta Appl. Math.
{\bf 44} (1996), pp. 295--352.
\bibitem[13]{Koe2} E. Koelink, {\it Eight lectures on quantum groups
and $q$-special functions}, Rev. Colombiana Mat. {\bf 30} (1996), no. 2,
pp. 93--180. 
\bibitem[14]{KR} E. Koelink, H. Rosengren, {\it Harmonic analysis on the 
$\hbox{SU}(2)$ dynamical quantum group}, Acta Appl. Math. {\bf 69}
(2001), pp. 163--220.
\bibitem[15]{KooIndag} T.H. Koornwinder, {\it Representations of the twisted
$\hbox{SU}(2)$ quantum group and some $q$-hypergeometric orthogonal 
polynomials}, Proc. Kon. Ned. Akad. van Wetensch., Ser. A {\bf 92}
(Indag. Math. {\bf 51}) (1989), pp. 97--117.
\bibitem[16]{K} T.H. Koornwinder, {\it Askey-Wilson polynomials
for root systems of type $BC$}, in ``Hypergeometric functions on domains
of positivity, Jack polynomials, and applications'' (Tampa, FL, 1991),
pp. 189--204, Contemp Math., {\bf 138}, Amer. Math. Soc., Providence, RI, 1992.
\bibitem[17]{Koo1} T.H. Koornwinder, {\it Askey-Wilson polynomials
as zonal spherical functions on the $\hbox{SU}(2)$ quantum group},
SIAM J. Math. Anal. {\bf 24} (1993), no. 3, pp. 795--813.
\bibitem[18]{Koo2} T.H. Koornwinder, {\it Some details of proofs of
theorems related to the quantum dynamical Yang-Baxter equation},
Int. J. Math. Math. Sci. {\bf 24} (2000), no. 12, pp. 793--806.
\bibitem[19]{KT} T.H. Koornwinder, N. Touhami, {\it Fusion and
exchange matrices for quantized $\hbox{sl}(2)$ and associated
$q$-special functions}, preprint (2002), math.QA/0207159.
\bibitem[20]{Mac} I.G. Macdonald, {\it Orthogonal polynomials associated
with root systems}, S{\'e}m. Lothar. Combin. {\bf 45} (2000/01),
Art. B45a.
\bibitem[21]{M} T. Masuda, K. Mimachi, Y. Nakagami, M. Noumi,
K. Ueno, {\it Representations of the quantum group $\hbox{SU}_q(2)$
and the little $q$-Jacobi polynomials}, J. Funct. Anal. {\bf 99} (1991),
pp. 357--386. 
\bibitem[22]{N} M. Noumi, {\it Macdonald's symmetric polynomials as zonal
spherical functions on some quantum homogeneous spaces}, Adv. Math.
{\bf 123} (1996), no. 1, pp. 16--77.
\bibitem[23]{NDS} M. Noumi, M.S. Dijkhuizen, T. Sugitani, {\it
Multivariable Askey-Wilson polynomials and quantum complex Grassmannians},
in ``Special functions, $q$-series and related topics''
(Toronto, ON, 1995), pp. 167--177, Fields Inst. Commun., {\bf 14},
Amer. Math. Soc., Providence, RI, 1997.
\bibitem[24]{NM} M. Noumi, K. Mimachi, {\it Askey-Wilson polynomials
and the quantum group $\hbox{SU}_q(2)$}, Proc. Japan Acad., Ser. A
{\bf 66} (1990), pp. 146--149.
\bibitem[25]{R0} H. Rosengren, {\it A new quantum algebraic interpretation
of the Askey-Wilson polynomials}, Contemp. Math. {\bf 254}
(2000), pp. 371--394.
\bibitem[26]{VS} L.L. Vaksman, Ya. S. Soibelman, {\it
Algebra of functions on the quantum group $\hbox{SU}(2)$},
Funct. Anal. Appl. {\bf 22} (1988), pp. 170--181.
\end{thebibliography}
\end{document}